\documentclass[review]{siamart190516}   
\usepackage{bm}

\usepackage{siampaper}



\newcommand{\cahiernumber}{17}  

\title{%
  An Interior-Point Trust-Region Method for Nonsmooth Regularized Bound-Constrained Optimization
}

\author{
  Geoffroy Leconte
  \thanks{%
    GERAD and Department of Mathematics and Industrial Engineering,
    Polytechnique Montr\'eal, QC, Canada.
    E-mail: \mailto{geoffroy.leconte@polymtl.ca}.
  }
  \and
  Dominique Orban%
  \thanks{%
    GERAD and Department of Mathematics and Industrial Engineering,
    Polytechnique Montr\'eal, QC, Canada.
    E-mail: \mailto{dominique.orban@gerad.ca}.
    Research partially supported by an NSERC Discovery Grant.
  }
}
\date{\today}

\usepackage{graphicx}
\graphicspath{{figs/}}

\usepackage{xcolor}
\usepackage{tikz}
\usepackage{tikzscale}
\usetikzlibrary{external}
\newcommand*{\includetikzgraphics}[2][]{%
	\includegraphics[#1]{#2}
}

\makeatletter
\renewcommand{\todo}[2][]{\tikzexternaldisable\@todo[#1]{#2}\tikzexternalenable}
\makeatother

\usepackage{pgfplots}
\pgfplotsset{compat=newest}
\pgfplotsset{
  tick label style={font=\scriptsize},
  title style={font=\scriptsize},
  label style={font=\scriptsize},
  legend style={font=\scriptsize},
  height = {0.48\linewidth},
  width = {0.48\linewidth},
  xticklabel style = {font=\scriptsize},
  yticklabel style = {font=\scriptsize},
  table/search path={figs},
}

\newtheorem{modelassumption}{Model Assumption}[section]
\newtheorem{stepassumption}{Step Assumption}[section]
\newtheorem{parameterassumption}{Parameter Assumption}[section]

\crefname{problemassumption}{Problem Assumption}{Problem Assumptions}
\Crefname{problemassumption}{Problem Assumption}{Problem Assumptions}
\crefname{modelassumption}{Model Assumption}{Model Assumptions}
\Crefname{modelassumption}{Model Assumption}{Model Assumptions}
\crefname{stepassumption}{Step Assumption}{Step Assumptions}
\Crefname{stepassumption}{Step Assumption}{Step Assumptions}
\crefname{parameterassumption}{Parameter Assumption}{Parameter Assumptions}
\Crefname{parameterassumption}{Parameter Assumption}{Parameter Assumptions}

\Crefname{subsection}{Section}{Sections}
\crefname{subsection}{section}{sections}

\newcommand{\dist}{\mathop{\textup{dist}}}
\newcommand{\diag}{\mathop{\textup{diag}}}

\newcommand{\cl}{\mathop{\textup{cl}}}
\newcommand{\epi}{\mathop{\textup{epi}}}
\newcommand{\gph}{\mathop{\textup{gph}}}
\newcommand{\lev}{\mathop{\textup{lev}}}
\newcommand{\eliminf}{\mathop{\textup{e-lim inf}}}
\newcommand{\elimsup}{\mathop{\textup{e-lim sup}}}
\newcommand{\elim}{\mathop{\textup{e-lim}}}
\newcommand{\glimsup}{\mathop{\textup{g-lim sup}}}
\newcommand{\gliminf}{\mathop{\textup{g-lim inf}}}
\newcommand{\glim}{\mathop{\textup{g-lim}}}
\newcommand{\pliminf}{\mathop{\textup{p-lim inf}}}
\newcommand{\plimsup}{\mathop{\textup{p-lim sup}}}
\newcommand{\plim}{\mathop{\textup{p-lim}}}

\newcommand{\B}{\mathds{B}}

\newcommand{\pto}{\overset{\textup{p}}{\longrightarrow}}
\newcommand{\eto}{\overset{\textup{e}}{\longrightarrow}}
\newcommand{\gto}{\overset{\textup{g}}{\longrightarrow}}

\makeatletter
\newcommand{\pushright}[1]{\ifmeasuring@#1\else\omit\hfill$\displaystyle#1$\fi\ignorespaces}

\begin{document}

  \nolinenumbers
  \maketitle

  \thispagestyle{firstpage}
  \pagestyle{myheadings}

  \begin{abstract}
    We develop an interior-point method for nonsmooth regularized bound-constrained optimization problems.
    Our method consists of iteratively solving a sequence of unconstrained nonsmooth barrier subproblems.
    We use a variant of the proximal quasi-Newton trust-region algorithm TR of \citet{aravkin-baraldi-orban-2022} to solve the barrier subproblems, with additional assumptions inspired from well-known smooth interior-point trust-region methods.
    We show global convergence of our algorithm with respect to the criticality measure of \citet{aravkin-baraldi-orban-2022}.
    Under an additional assumption linked to the convexity of the nonsmooth term in the objective, we present an alternative interior-point algorithm with a slightly modified criticality measure, which performs better in practice.
    Numerical experiments show that our algorithm performs better than the trust-region method TR, the trust-region method with diagonal hessian approximations TRDH of \citet{leconte-orban-2023}, and the quadratic regularization method R2 of \citet{aravkin-baraldi-orban-2022} for two out of four tested bound-constrained problems. 
    On those two problems, our algorithm obtains smaller objective values than the other solvers using fewer objective and gradient evaluations.
    On the two other problems, it performs similarly to TR, R2 and TRDH.
  \end{abstract}



\section{Introduction}%
\label{sec:intro}

We consider the problem
\begin{equation}%
  \label{eq:nlp}
  \minimize{x \in \R^n} f(x) + h(x) \quad \st \ x \geq 0,
\end{equation}
where \(f: \R^n \to \R\) has Lipschitz-continuous gradient with constant \(L_f \ge 0\) and \(h: \R^n \to \R \cup \{+\infty\}\) is proper and lower semi-continuous.
Both \(f\) and \(h\) may be nonconvex.
\(h\) is often considered as a regularization function used to favor solutions with desirable properties, such as sparsity.

Problems such as~\eqref{eq:nlp} are classically solved with a variant of the proximal gradient method \citep{lions-mercier-1979}.
The proximal quasi-Newton trust-region algorithm of \citet{aravkin-baraldi-orban-2022}, referred to as TR, can be extended to box constraints, provided that the proximal operator of \(h + \chi(\cdot \mid [\ell, u])\), where \(\ell < u\) componentwise and \(\chi\) is the indicator of the box \([\ell, u]\), can be computed efficiently.
\citet{leconte-orban-2023} present a variant of TR named TRDH that also supports box constraints, and uses diagonal quasi-Newton approximations.
For specific separable regularizers \(h\), they provide a closed-form solution of the trust-region subproblems with box constraints, giving rise to what they defined as an indefinite proximal operator; a scaled generalization of a proximal operator.

At each iteration \(k\), we solve the barrier subproblem
\begin{equation}
  \label{eq:subprob-intro}
  \minimize{x \in \R^n} f(x) + \phi_k(x) + h(x),
\end{equation}
where \(\phi_k\) is the logarithmic barrier function 
\begin{equation}
  \label{eq:log-barrier-func}
  \phi_k(x) = \sum_{i=1}^n \phi_{k,i}(x), \qquad \phi_{k,i}(x) := -\mu_k \log(x_i), \ i = 1, \ldots, n,
\end{equation}
and \(\{\mu_k\} \searrow 0\).
Each subproblem~\eqref{eq:subprob-intro} is an unconstrained problem if we consider that \(-\log(x) = +\infty\) when \(x \le 0\).
In \Cref{sec:convergence-inner}, we explain that under reasonable assumptions, we can solve solve~\eqref{eq:subprob-intro} with a modified version of \citeauthor{aravkin-baraldi-orban-2022}'s TR algorithm, and we expect that the solutions of~\eqref{eq:log-barrier-func} converge to a solution of~\eqref{eq:nlp} as \(\mu_k \to 0\).

Our approach is sometimes referred to as a trust-region interior-point method, or trust-region method for barrier functions. 
We refer the reader to \citep[Chapter~\(13\)]{conn-gould-toint-2000} for more information on the case where \(h = 0\).
Our algorithm, named RIPM (\emph{Regularized Interior Proximal Method}), can be seen as a generalization of those methods to solve~\eqref{eq:nlp}.

An inconvenient of solving~\eqref{eq:subprob-intro}, induced by the logarithmic barrier function~\eqref{eq:log-barrier-func}, is that the smooth part of the subproblem \(f + \phi_k\) does not have a Lipschitz gradient, thus compromising the convergence properties of TR established by \citet{aravkin-baraldi-orban-2022}.
Nevertheless, in our analysis, we establish the convergence of the barrier subproblems using the update rules of \citet[Chapter~\(13.6.3\)]{conn-gould-toint-2000} for our trust-region model.
We also show global convergence of RIPM towards a first-order stationary point of~\eqref{eq:nlp} if the trust-region radii and the step lengths used in proximal operator evaluations are bounded away from zero, and the iterates generated by the algorithm remain bounded.
In \Cref{sec:new-crit-meas}, under a convexity assumption on the nonsmooth term, we provide an alternative implementation of the outer iterations of RIPM where we change the stopping criteria to improve numerical performance.

In addition, we implement a variant of RIPM named RIPMDH (\emph{Regularized Interior Proximal Method with Diagonal Hessian approximations}) that uses TRDH to solve the barrier subproblems.
We compare the performance of RIPM and RIPMDH with TR, TRDH and R2, all available from \href{https://github.com/JuliaSmoothOptimizers/RegularizedOptimization.jl}{RegularizedOptimization.jl} \citep{baraldi-orban-regularized-optimization-2022}, on four bound-constrained problems.
The first two problems are a regularized box-constrained quadratic problem, and a sparse nonnegative matrix factorization problem.
These two problems require many TR and R2 iterations to converge.
RIPM and RIPMDH obtain smaller objective values than the other solvers using fewer objective and gradient evaluations, which suggests that they may be best suited to solve difficult bound-constrained nonsmooth problems.
The third problem is an inverse problem for finding the parameters of a differential equation.
RIPM and RIPMDH perform more objective and gradient evaluations than TR, but RIPMDH performs the least amount of proximal operator evaluations.
The last problem is a regularized box-constrained basis-pursuite denoise problem.
RIPMDH exhibits similar performance to those of TR, TRDH and R2 using a modification of some of its parameters.

\subsubsection*{Related research}

\citet{attouch-wets-1981} use \(f = 0\) and the nonsmooth barrier \(-\mu_k \sum_{i=1}^n \log(\min(\tfrac{1}{2}, x_i)) > 0\).
They use the theory of epi-convergence to explain some convergence properties of barrier methods, and in particular, that the objectives of the barrier subproblems epi-converge to the objective of their initial constrained problem.

\citet{chouzenoux-corbineau-pesquet-2020} use convex \(f\) and \(h\) with general inequality constraints \(c_i(x) \le 0\) for \(i \in \{1, ..., p\}\), \(p > 0\) to solve large-scale image processing problems.
Their algorithm uses proximal gradient steps to solve the barrier subproblem.

\citet{bertocchi-chouzenoux-corbineau-pesquet-prato-2020} solve inverse problems \(y = \mathcal{D}(H \bar x)\), where \(y\) is some observed data, \(\bar x\) the signal to determine, \(H\) is a linear observation operator, and \(\mathcal{D}\) is an operator applying noise perturbation.
They solve the problem \(\minimize{x \in \R^n} f(Hx, y) + h(x)\) with constraints \(\ c_i(x) \geq 0, i \in \{1, ..., p\}\), where \(f\) is a convex function involving the observations and the signal \(x\), \(f(\cdot; y)\) and \(h\) are twice differentiable, \(h\) and \(-c_i\) are convex (among other properties).
They compute the proximal operator of the barrier term, and use an interior-point algorithm.
They apply their algorithm to develop a neural network architecture for image restoration.  

\citet{demarchi-themelis-2022} use the proximal gradient method to solve nonsmooth regularized optimization problems such as~\eqref{eq:nlp} where \(f\) has a locally Lipschitz-continuous gradient, \(h\) is continuous relative to its domain and prox-bounded.
In addition, the constraint \(x \ge 0\) is replaced by a more general constraint \(c(x) \le 0\), where \(c\) has locally Lipschitz-continuous Jacobian.

\citet{shen-xue-zhang-wang-2020} present an active set proximal algorithm to solve~\eqref{eq:nlp} with \(h(x) = \lambda \|x\|_1\) for some \(\lambda > 0\), and where \(\ell \le x \le u\) with \(\ell < 0 < u\) instead of \(x \ge 0\).
They use a hybrid search direction based upon a proximal gradient step for the active variables (i.e., the variables that are at one the bounds of the constraints), and a Newton step for the other variables.

\subsubsection*{Notation}

For \(v \in \R^n\), \(\|v\|\) denotes the Euclidean norm of \(v\).
\(\R_+\) and \(\R_{++}\), denote, respectively, the sets of positive and strictly positive real numbers, whereas \(\R_+^n\) and \(\R_{++}^n\) denote the sets of vectors having all their components in \(\R_+\) and \(\R_{++}\), respectively.

\(\widebar \R\) denotes \(\R \cup \{\pm \infty\}\).
The unit closed ball defined with the \(\ell_{\infty}\)-norm and centered at the origin is \(\B\), and the ball centered at the origin of radius \(\Delta > 0\) is \(\Delta \B\).
If \(C \subseteq \R^n\), the indicator of \(C\) is \(\chi(\cdot \mid C): \R^n \to \widebar{\R}\) defined by \(\chi(x \mid C) = 0\) if \(x \in C\), and \(\chi(x \mid C) = +\infty\) if \(x \not \in C\).
For \(y \in \R^n\), the set \(y + C\) is composed of all the vectors \(s \in \R^n\) such that \(s = y + x\) with \(x \in C\). 

Following the notation of \citet{rtrw}, the set of all subsequences of \(\N\) is denoted by \(\mathcal{N}_{\infty}^{\#}\), and the set composed of the subsequences of \(\N\) containing all \(k\) beyond some \(k_0\) is denoted by \(\mathcal{N}_{\infty}\).
For \(N \in \mathcal{N}_{\infty}^{\#}\), \(\{x_k\} \underset{N}{\longrightarrow} \bar x\) indicates that the subsequence \(\{x_k\}_{k \in N}\) (which we may also write \(\{x_k\}_N\) for conciseness) converges to \(\bar x\).

\(X\), \(Z\) and \(S\) (possibly with subscripts\(~_k\),\(~_{k,j}\) or\(~_{k,j,1}\)) denote the square diagonal matrices having \(x\), \(z\) and \(s\) as diagonal elements, respectively.

\section{Background}%
\label{sec:background}

The following are standard variational analysis concepts---see, e.g., \citep{rtrw}.
Let \(\phi: \R^n \to \widebar{\R}\) and \(\bar{x} \in \R^n\) where \(\phi\) is finite.
The Fr\'echet subdifferential of \(\phi\) at \(\bar{x}\) is the closed convex set \(\widehat{\partial} \phi(\bar{x})\) of elements \(v \in \R^n\) such that
\[
  \liminf_{\substack{x \to \bar{x} \\ x \neq \bar{x}}} \frac{\phi(x) - \phi(\bar{x}) - v^T (x - \bar{x})}{\|x - \bar{x}\|} \geq 0.
\]
The limiting subdifferential of \(\phi\) at \(\bar{x}\) is the closed, but not necessarily convex, set \(\partial \phi(\bar{x})\) of elements \(v \in \R^n\) for which there exist \(\{x_k\} \to \bar{x}\) and \(\{v_k\} \to v\) such that \(\{\phi(x_k)\} \to \phi(\bar{x})\) and \(v_k \in \widehat{\partial} \phi(x_k)\) for all \(k\).
The inclusion \(\widehat{\partial} \phi(\bar{x}) \subset \partial \phi(\bar{x})\) always holds.
Finally, the horizon subdifferential of \(\phi\) at \(\bar{x}\) is the closed, but not necessarily convex, cone \(\partial^{\infty} \phi(\bar{x})\) of elements \(v \in \R^n\) for which there exist \(\{x_k\} \to \bar{x}\), \(\{v_k\}\) and \(\{\lambda_k\} \searrow 0\) such that \(\{\phi(x_k)\} \to \phi(\bar{x})\), \(v_k \in \widehat{\partial} \phi(x_k)\) for all \(k\), and \(\{\lambda_k v_k\} \to v\).

If \(\inf \phi > -\infty\), \(\argmin{} \phi\) is the set of \(x \in \R^n\) such that \(\phi(x) = \inf \phi\).
For \(\epsilon > 0\), \(\epsilon\text{-}\argmin{} \phi\) is the set of \(x \in \R^n\) such that \(\phi(x) \leq \inf \phi + \epsilon\).

If \(C \subseteq \R^n\) and \(\bar{x} \in C\), the closed convex cone \(\widehat{N}_C(\bar{x}) := \widehat{\partial} \chi(\bar{x} \mid C)\) is the regular normal cone to \(C\) at \(\bar{x}\).
The closed cone \(N_C(\bar{x}) := \partial \chi(\bar{x} \mid c) = \partial^{\infty} \chi(\bar{x} \mid C)\) is the normal cone to \(C\) at \(\bar{x}\).
\(\widehat{N}_C(\bar{x}) \subseteq N_C(\bar{x})\) always holds, and is an equality if \(C\) is convex.

If \(C\) is convex, \(N_C(\bar{x})\) is the closed convex cone of elements \(v \in \R^n\) such that \(v^T (x - \bar{x}) \leq 0\) for all \(x \in C\) \citep[Theorem~\(6.9\)]{rtrw}.

For a set-valued mapping \(S: X \rightarrow U\) where, for any \(x \in X\), \(S(x) \subset U\), the graph of \(S\) is the set \(\gph S := \{(x, u) \mid u \in S(x)\}\).

For \(\bar x \in \R^n\), the limit superior of \(S\) at \(\bar x\) is \(\limsup_{x \rightarrow \bar x} S(x) := \{u \mid \exists \ \{x_k\} \rightarrow \bar x, \exists \ \{u_k\} \rightarrow u \text{ with } u_k \in S(x_k)\}\), and the limit inferior of \(S\) at \(\bar x\) is \(\liminf_{x \rightarrow \bar x} S(x) := \{u \mid \forall \ \{x_k\} \rightarrow \bar x, \exists \ N \in \mathcal{N}_{\infty}, \{u_k\} \underset{N}{\rightarrow} u \text{ with } u_k \in S(x_k)\}\).
\(S(\bar x) \subseteq \limsup_{x \rightarrow \bar x} S(x)\) and \(S(\bar x) \supseteq \liminf_{x \rightarrow \bar x} S(x)\) always hold.

The set-valued mapping \(S\) is outer semicontinuous (osc) at \(\bar x\) if \(\limsup_{x \rightarrow \bar x} S(x) \subseteq S(\bar x)\), or, equivalently, \(\limsup_{x \rightarrow \bar x} S(x) = S(\bar x)\).
It is inner semicontinuous (isc) at \(\bar x\) if \(\liminf_{x \rightarrow \bar x} S(x) \supseteq S(\bar x)\), or equivalently \(\liminf_{x \rightarrow \bar x} S(x) = S(\bar x)\) when \(S\) is closed-valued.
If both conditions hold, \(S\) is continuous at \(\bar x\), i.e., \(S(x) \rightarrow S(\bar x)\) as \(x \rightarrow \bar x\). 

\begin{shadyproposition}[{\protect \citealp[Proposition~\(8.7\)]{rtrw}}]
  \label{prop:subgradient-osc}
  For \(\phi: \R^n \to \widebar{\R}\) and \(\bar x\) where \(\phi\) is finite, \(\partial \phi\) is osc at \(\bar x\) with respect to \(\phi(x) \rightarrow \phi(\bar x)\) when \(x \rightarrow \bar x\), i.e. for any \(\{x_k\} \rightarrow \bar x\) with \(\{\phi(x_k)\} \rightarrow \phi (\bar x)\), there exists \(v_k \in \partial \phi (x_k)\) for all \(k\) such that \(\{v_k\} \rightarrow \bar v \in \partial \phi(\bar x)\).
\end{shadyproposition}

The graphical outer limit of a sequence of set-valued mappings \(S_k\) is defined by \((\glimsup_k S_k)(x) := \{u \mid \exists N \in \mathcal{N}_{\infty}^{\#}, \{x_k\} \underset{N}{\rightarrow} x, \{u_k\} \underset{N}{\rightarrow} u, u_k \in S_k(x_k)\}\).
The graphical inner limit of a sequence of set-valued mappings \(S_k\) is defined by \((\gliminf_k S_k)(x) := \{u \mid \exists N \in \mathcal{N}_{\infty}, \{x_k\} \underset{N}{\rightarrow} x, \{u_k\} \underset{N}{\rightarrow} u, u_k \in S_k(x_k)\}\).
If both limits agree, the graphical limit \(S = \glim S_k\) exists, so that we can also write \(S_k \gto S\), and we have \(S_k \gto S \iff \gph S_k \rightarrow \gph S\).

The epigraph of \(\phi\) is the set \(\epi \phi = \{ (x, \alpha) \in \R^n \times \R \mid \alpha \geq \phi(x) \}\).

We denote \(\cl(\phi)\) the (lower) closure of \(\phi\), i.e., the largest function less than \(\phi\) that is lower semi-continuous.
Its epigraph is the closure of \(\epi \phi\).

If \(\phi_k: \R^n \to \widebar{\R}\) for all \(k \geq 0\), the lower and upper pointwise limits of \(\{\phi_k\}\) are the functions \(\pliminf_k \phi_k\) and \(\plimsup_k \phi_k: \R^n \to \widebar{\R}\) defined for all \(x \in \R^n\) by
\begin{align*}
  (\pliminf_k \phi_k)(x) & := \liminf_k \phi_k(x),
  \\ (\plimsup_k \phi_k)(x) & := \limsup_k \phi_k(x).
\end{align*}
When \(\pliminf_k \phi_k\) and \(\plimsup_k \phi_k\) coincide, their common value is the pointwise limit of \(\{\phi_k\}\) denoted \(\plim_k \phi_k\).
If \(\plim_k \phi_k = \phi\), we write \(\{\phi_k\} \pto \phi\).

For a sequence \(\{E_k\} \subseteq \R^n\), we define
\begin{align*}
  \limsup_k E_k & := \{z \in \R^n \mid \exists N \in \mathcal{N}_{\infty}^{\#}, \exists \{z_j\}_N \to z, \ z_j \in E_j \text{ for all } j \in N \}
  \\ \liminf_k E_k & := \{z \in \R^n \mid \exists N \in \mathcal{N}_{\infty}, \exists \{z_j\}_N \to z, \ z_j \in E_j \text{ for all } j \in N \}.
\end{align*}

Consider in particular \(E_k := \epi \phi_k\).
It is not difficult to see that \(\limsup_k E_k\) and \(\liminf_k E_k\) are also epigraphs in the sense that if \((x, t)\) is in either set, then so is \((x, s)\) for any \(s \geq t\).
The lower epi-limit of \(\{\phi_k\}\) is the function \(\eliminf_k \phi_k\) whose epigraph is \(\limsup_k E_k\), and the upper epi-limit of \(\{\phi_k\}\) is the function \(\elimsup_k \phi_k\) whose epigraph is \(\liminf_k E_k\).
It is always true that \(\eliminf_k \phi_k \leq \elimsup_k \phi_k\).
When the two coincide, their common value is called the epi-limit of \(\{\phi_k\}\) denoted \(\elim_k \phi_k\).
If \(\elim_k \phi_k = \phi\), we also write \(\{\phi_k\} \eto \phi\).

The following result summarizes important properties of epi-limits used in the sequel.

\begin{shadyproposition}[{\protect \citealp[Proposition~\(7.4\)]{rtrw}}]
  Let \(\phi_k: \R^n \to \R\) for \(k \geq 0\).
  \begin{enumerate}
    \item \(\eliminf_k \phi_k\) and \(\elimsup_k \phi_k\) are lsc, and so is \(\elim_k \phi_k\) when it exists;
    \item if \(\phi_k \geq \phi_{k+1}\) for all \(k\), \(\elim_k \phi_k\) exists and equals \(\cl (\inf_k \phi_k)\);
    \item if \(\phi_k \leq \phi_{k+1}\) for al \(k\), \(\elim_k \phi_k\) exists and equals \(\cl (\sup_k \phi_k)\).
  \end{enumerate}
  In addition, if \(\phi\), \(\underline{\phi}_k\), \(\widebar{\phi}_k: \R^n \to \R\) with \(\underline{\phi}_k \leq \phi_k \leq \widebar{\phi}_k\) for all \(k\), and if \(\{\underline{\phi}_k\} \eto \phi\) and \(\{\widebar{\phi}_k\} \eto \phi\), then \(\{\phi_k\} \eto \phi\).
\end{shadyproposition}

The model that we will use in our algorithm uses an approximation \(\psi (\cdot; x)\) of \(h\) at \(x\) so that \(\psi(s, x) \approx h(x + s)\).
For \(\psi : \R^n \times \R^m \rightarrow \widebar \R\), the function-valued mapping \(x \mapsto \psi(\cdot; x)\) is \emph{epi-continous} at \(\bar x\) if \(\psi(\cdot, x) \eto \psi(\cdot, \bar x)\) as \(x \rightarrow \bar x\).

\(\phi\) is level-bounded if, for every \(\alpha \in \R\), the lower level set \(\lev_{\le \alpha} \phi := \{x \in \R^n \mid \phi(x) \le \alpha\}\) is bounded (possibly empty).
The sequence of functions \(\{\phi_k\}\) is eventually level-bounded if, for each \(\alpha \in \R\), the sequence of sets \(\{\lev_{\le \alpha} \phi_k\}\) is eventually bounded, i.e., there is an index set \(N \in \mathcal{N}_{\infty}\) such that \(\{\lev_{\le \alpha} \phi_k\}_{k \in N}\) is bounded.

The following theorem establishes properties about the minimization of sequences of epi-convergent functions. 

\begin{shadytheorem}[{\protect \citealp[Theorem~\(7.33\)]{rtrw}}]
  \label{thm:epi-cv-minimization}
  Suppose the sequence \(\{\phi_k\}\) is eventually level-bounded, and \(\phi_k \eto \phi\) with \(\phi_k\) and \(\phi\) lsc and proper.
  Then,
  \begin{equation}
    \label{eq:cv-inf-phik}
    \inf \phi_k \rightarrow \inf \phi
  \end{equation}
  with \(-\infty < \inf \phi < +\infty\), while there exists \(N \in \mathcal{N}_{\infty}\) such that \(\argmin{} \phi_k\) is a bounded sequence of nonempty sets with
  \begin{equation}
    \label{eq:limsup-argmin-phik}
    \limsup_k (\argmin{} \phi_k) \subset \argmin{} \phi.
  \end{equation}
  Indeed, for any \(\{\epsilon_k\} \searrow 0\) and \(x_k \in \epsilon_k\text{-}\argmin{} \phi_k\), \(\{x_k\}\) is bounded and all its cluster points belong to \(\argmin{} \phi\).
  If \(\argmin{} \phi\) consists of a unique point \(\bar x\), one must actually have \(\{x_k\} \rightarrow \bar x\).
\end{shadytheorem}

The proximal operator associated with the proper lsc function \(h\) and parameter \(\nu > 0\) is 
\begin{equation}
  \prox{\nu h}(x) := \arg \min_w \tfrac{1}{2}\nu^{-1}\|w - x\|^2 + h(w).
\end{equation}

  \section{Stationarity}%
\label{sec:stationarity}

First-order stationarity conditions for~\eqref{eq:nlp} may be stated as \citep[Theorem~\(10.1\)]{rtrw}
\begin{equation}%
  \label{eq:nlp-stationarity}
  0 \in \nabla f(x) + \partial (h + \chi(\cdot \mid \R^n_+))(x).
\end{equation}

We say that the constraint qualification (CQ) holds at \(x\) if \(\partial^{\infty} (f + h)(x)\) contains no \(v \neq 0\) such that \(-v \in N_{\R^n_+}(x)\).

Under the constraint qualification,~\eqref{eq:nlp-stationarity} can also be written \citep[Theorem~\(8.15\)]{rtrw}
\begin{equation}
  \label{eq:nlp-stationarity-bis}
  0 \in \nabla f(x) + \partial h(x) + N_{\R^n_+}(x).
\end{equation}

Our assumptions that \(f\) is continuously differentiable and that \(h\) is proper lsc allows us to write \citep[Exercise~\(8.8\) and Theorem~\(8.9\)]{rtrw}
\[
  \partial^{\infty}(f + h)(x) =
  \partial^{\infty} h(x) =
  \{v \in \R^n \mid (v, 0) \in N_{\epi h}(x, h(x)) \}.
\]
Due to the simple form of \(N_{\R^n_+}\), the constraint qualification states that the only \(v \in \partial^{\infty} h(x)\) such that \(v \geq 0\) and satisfying \(v_i = 0\) if \(x_i > 0\) is \(v = 0\).

A simple example where the constraint qualification is not satisfied is found by setting \(n = 1\), \(f = 0\), and \(h(x) = |x|_0\), i.e., \(h(x) = 1\) if \(x \neq 0\) and \(h(0) = 0\), in~\eqref{eq:nlp}.
The unique solution is \(\bar{x} = 0\).
Then, \(\widehat{N}_{\epi h}(0, 0) = N_{\epi h}(0, 0) = \{(v, t) \mid t \leq 0\}\) so that \(\partial^{\infty} h(0) = \R\).
The qualification condition requires that the only \(v \in \R\) such that \(v \geq 0\) be \(v = 0\), which is clearly not the case.
Of course, the bound constraint in the example above is redundant and the constrained and unconstrained solutions coincide.

Using \citep[Example~\(6.10\)]{rtrw}, \(N_{\R_+^n}(x) = N_{\R_+}(x_1) \times \dots \times N_{\R_+}(x_n)\), where \(N_{\R_+}(0) = (-\infty, 0]\) and for all \(x_i > 0\), \(N_{\R_+}(x_i) = \{0\}\).
Thus,~\eqref{eq:nlp-stationarity-bis} can also be formulated as
\begin{equation}
  \label{eq:nlp-stationarity-z}
  0 \in \nabla f(x) + \partial h(x) - z, \quad Xz = 0, \quad x \ge 0, \quad z \ge 0,
\end{equation}
where \(X = \diag(x)\) and \(Z = \diag(z)\).

For fixed \(x \in \R^n_+\) and \(z \in \R_+^n\), we define approximations
\begin{subequations}
  \label{eq:def-modelbar}
  \begin{align}
    \varphi^{\mathcal{L}}(s; x, z) & := f(x) + (\nabla f(x) - z)^T s, \label{eq:bar-varphi} \\
    \psi(s; x) & \phantom{:}\approx h(x + s) \text{ with } \psi(0; x) = h(x) \text{ and } \partial \psi(0; x) = \partial h(x),
    \label{eq:def-psi} \\
    \hat{\psi}(s; x) & := \psi(s; x) + \chi(x + s \mid \R^n_+),
    \label{eq:def-psibar}
  \end{align}
\end{subequations}
and the model of \(f + h\) about \(x\)
\begin{equation}
    m^{\mathcal{L}}(s; x, z, \nu) := \varphi^{\mathcal{L}}(s; x, z) + \tfrac{1}{2}\nu^{-1}\|s\|^2 + \psi(s; x), \label{eq:bar-model}
\end{equation}
where \(\nu > 0\).
We point out that \(\nabla \varphi^{\mathcal{L}}(s; x, z) = \nabla f(x) - z\), which is the expression of the Lagrangian in the smooth case, thus, we use the superscript\(~^{\mathcal{L}}\) to denote objects sharing similarities with the smooth Lagrangian.

For \(\Delta \geq 0\), we further define
\begin{subequations}
  \begin{align}
    p^{\mathcal{L}}(\Delta; x, z, \nu) & := \min_s \ \varphi^{\mathcal{L}}(s; x, z) + \tfrac{1}{2}\nu^{-1}\|s\|^2 + \hat{\psi}(s; x) + \chi(s \mid \Delta \B),
    \label{eq:def-pbar} \\
    P^{\mathcal{L}}(\Delta; x, z, \nu) & := \argmin{s} \ \varphi^{\mathcal{L}}(s; x, z) + \tfrac{1}{2}\nu^{-1}\|s\|^2 + \hat{\psi}(s; x) + \chi(s \mid \Delta \B).
    \label{eq:def-Pbar}
  \end{align}
\end{subequations}
Our associated optimality measure is
\begin{equation}
  \label{eq:def-xibar}
  \xi^{\mathcal{L}}(\Delta; x, z, \nu) := f(x) + h(x) - \varphi^{\mathcal{L}}(s^{\mathcal{L}}; x, z) - \psi(s^{\mathcal{L}}; x),
\end{equation}
where \(s^{\mathcal{L}} \in P^{\mathcal{L}}(\Delta, x, z, \nu)\).

\begin{shadylemma}
  \label{lem:xi-stationarity}
  Let the CQ hold at \(x \in \R^n_+\), \(z \ge 0\) such that \(Xz = 0\) and \(\Delta > 0\).
  Then, \(\xi^{\mathcal{L}}(\Delta; x, z, \nu) = 0\) \(\Longleftrightarrow\) \(0 \in P^{\mathcal{L}}(\Delta; x, z, \nu)\) \(\Longrightarrow\) \(x\) is first-order stationary for~\eqref{eq:nlp}.
\end{shadylemma}

\begin{proof}
  The first equivalence follows directly from~\eqref{eq:def-Pbar}--\eqref{eq:def-xibar}.
  The first-order necessary conditions for~\eqref{eq:def-pbar} then imply
  \begin{equation}
    \label{eq:first-order-cond-s0}
    \begin{aligned}
      0 & \in \nabla \varphi^{\mathcal{L}}(0; x, z) + \partial (\hat{\psi}(\cdot; x) + \chi(\cdot \mid \Delta \B))(0) \\
        & = \nabla f(x) - z + \partial (\psi(\cdot; x) + \chi(\cdot \mid (-x + \R_+^n)) + \chi(\cdot \mid \Delta \B))(0) \\
        & = \nabla f(x) - z + \partial (\psi(\cdot; x) + \chi(\cdot \mid (-x + \R_+^n) \cap \Delta \B))(0).
    \end{aligned}
\end{equation}
  As \((-x + \R_+^n)\) and \(\Delta \B\) are convex, so is \((-x + \R_+^n) \cap \Delta \B\).
  From this observation, we deduce that,
  \begin{align*}
    N_{(-x + \R^n_+) \cap \Delta \B}(0) & = \partial \chi(0 \mid (-x + \R^n_+) \cap \Delta \B) \\
                                        & = \partial \chi(0 \mid (-x + \R^n_+)) + \partial \chi(0 \mid \Delta \B) \\
                                        & = \partial \chi(x \mid \R^n_+) \\
                                        & = N_{\R_+^n}(x).
  \end{align*}
  The CQ combined with the above equations indicate that there is no \(v \in \partial^{\infty} h(x) = \partial^{\infty} \psi(0; x)\), \(v \neq 0\) such that \(-v \in N_{\R_+^n}(x) = N_{(-x + \R_+^n) \cap \Delta \B}(0)\), thus, \citep[Corollary~\(10.9\)]{rtrw} leads to
  \begin{equation}
    \label{eq:subgrad-addition}
    \begin{aligned}
      \partial (\psi(\cdot; x) + \chi(\cdot \mid (-x + \R_+^n) \cap \Delta \B))(0) & \subset \partial \psi(0; x) + N_{(-x + \R_+^n) \cap \Delta \B}(0) \\
                                                                                   & = \partial h(x) + N_{\R^n_+}(x).
    \end{aligned}
  \end{equation}
  By injecting~\eqref{eq:subgrad-addition} into~\eqref{eq:first-order-cond-s0}, we obtain
  \begin{equation*}
    0 \in \nabla f(x) - z + \partial h(x) + N_{\R^n_+}(x).
  \end{equation*}
  From the observation above~\eqref{eq:nlp-stationarity-z} and the fact that \(Xz = 0\), we deduce that for any \(v \in N_{\R^n_+}(x)\), \(v - z \in N_{\R_+^n}(x)\).
  Thus,
  \[
    0 \in \nabla f(x) + \partial h(x) + N_{\R_+^n}(x).
  \]
\end{proof}

If \(h\) is convex, the CQ is not required in \Cref{lem:xi-stationarity} \citep[Exercise~\(10.8\)]{rtrw}.

\section{Projected-directions methods}%
\label{sec:projdir}

Let us briefly recall the proximal gradient method \citep{lions-mercier-1979} used to solve
\begin{equation}
  \label{eq:nlp-unconstrained}
  \minimize{s \in \R^n} \ f(s) + \tilde h(s),
\end{equation}
where \(f: \R^n \to \R\) has Lipschitz-continuous gradient and \(\tilde h: \R^n \to \widebar{\R}\) is proper and lower semi-continuous.
The method generates iterates \(s_k\) such that
\begin{equation}
  \label{eq:pg-iter-general}
  s_{k+1} \in \prox{\nu \tilde h}(s_k - \nu \nabla f(s_k)) = \argmin{s} f(s_k) + \nabla f(s_k)^T (s - s_k) + \tfrac{1}{2} \nu^{-1}\|s - s_k\|^2 + \tilde h(s),
\end{equation}
where \(\nu > 0\), which leads to the first-order stationarity conditions
\begin{equation}
  \label{eq:prox-first-stat}
  0 \in s_{k+1} - s_k + \nu \nabla f(s_k) + \nu \partial \tilde h (s_{k+1}).
\end{equation}

A first approach to solving~\eqref{eq:nlp}, that we can reformulate as
\[
  \minimize{x \in \R^n} \ f(x) + h(x) + \chi(x \mid \R^n_+),
\]
is to use projected-directions methods.
A simple example of such methods consists in performing the identification \(\tilde h = h + \chi(\cdot \mid \R_+^n)\) in~\eqref{eq:nlp-unconstrained}, and using the proximal gradient method as in~\eqref{eq:pg-iter-general}.

\citeauthor{aravkin-baraldi-orban-2022}'s TR and R2 are other examples of algorithms that can solve~\eqref{eq:nlp} with a similar strategy.
Replacing \(h\) by \(\tilde h\) and \(\psi(\cdot; x)\) by \(\tilde \psi(\cdot; x) = \psi(\cdot; x) + \chi(\cdot \mid -x +\R_+^n)\) in all models of TR and R2 is sufficient to generalize these methods to~\eqref{eq:nlp}, if a solution of
\begin{equation}
  \label{eq:prox-proj}
  s_k \in \prox{\nu_k \tilde \psi(\cdot; x_k) + \chi(\cdot \mid  \Delta_k \B)}(-\nu_k \nabla f(x_k)).
\end{equation}
for TR, or 
\begin{equation}
  \label{eq:prox-proj-R2}
  s_k \in \prox{\nu_k \tilde \psi(\cdot; x_k)}(-\nu_k \nabla f(x_k))
\end{equation}
for R2, is available.

\citet{leconte-orban-2023} implement a variant of TR named TRDH that handles bound constraints for separable regularizers \(h\) (assuming that \(\psi\) is also separable).
TRDH solves at each iteration \(k\) and for all \(i \in \{1, \ldots, n\}\) the problem
\begin{equation}
  \label{eq:iprox-proj}
  (s_k)_i \in \arg \min_{s_i} \nabla f(x_k))_i s_i + \tfrac{1}{2} (d_k)_i s_i^2 + (\psi(s; x_k))_i + \chi(s_i \mid \Delta_k \B \cap (-(x_k)_i + \R_+),
\end{equation}
with \((d_k)_i \in \R\).
The special choice \((d_k)_i = \nu_k^{-1}\) shows that solving~\eqref{eq:iprox-proj} for all \(i\) is equivalent to solving~\eqref{eq:prox-proj}.

However, for nonseperable regularizers \(h\), projected-directions methods rely on computing search directions such as~\eqref{eq:prox-proj},~\eqref{eq:prox-proj-R2} or~\eqref{eq:iprox-proj}, which may be complicated (impossible for the latter), and therefore seems to be a limitation of this approach.
The following section describes the implementation of a method that is different from projected directions methods, and is based upon interior-point techniques.

\section{Barrier methods}%
\label{sec:barrier}

Consider a sequence \(\{\mu_k\} \searrow 0\).

\begin{shadylemma}%
  \label{lem:elim-bar}
  Let \(\phi_k\) be defined as in~\eqref{eq:log-barrier-func}.
  Then, \(\elim \phi_k = \chi(\cdot \mid \R^n_+)\).
\end{shadylemma}

\begin{proof}
  It is sufficient to show that \(\elim \phi_{k,i} = \chi(\cdot \mid \R_+)\) for \(i = 1, \ldots, n\).
  Our goal is to bound each \(\phi_k\) by two functions having \(\chi(\cdot \mid \R_+^n)\) as epi-limit.
  We define
  \[
    \phi_{k,i}^>(x) :=
    \begin{cases}
      +\infty & \text{if } x \leq 0 \\
      \phi_{k,i}(x) & \text{if } 0 < x < 1 \\
      0 & \text{if } x \geq 1,
    \end{cases}
    \qquad
    \phi_{k,i}^<(x) :=
    \begin{cases}
      +\infty & \text{if } x \leq 0 \\
      0  & \text{if } 0 < x < 1 \\
      \phi_{k,i}(x) & \text{if } x \geq 1.
    \end{cases}
  \]
  By construction, \(\phi_{k,i}(x) = \phi_{k,i}^>(x) + \phi_{k,i}^<(x)\), \(\{\phi_{k,i}^>(x)\} \searrow 0\) and \(\{\phi_{k,i}^<(x)\} \uparrow 0\) as \(k \to \infty\) for all \(x > 0\).
  In particular, \(\{\phi_{k,i}^>\} \pto \chi(\cdot \mid \R_+)\) and \(\{\phi_{k,i}^<\} \pto \chi(\cdot \mid \R_+)\) as \(k \to \infty\).

  By \citep[Proposition~\(7.4\)c]{rtrw}, because \(\{\phi_{k,i}^>\}\) is nonincreasing with \(k\), its epi-limit is well defined and \(\{\phi_{k,i}^>\} \eto \cl \inf_k \phi_{k,i}^> = \chi(\cdot \mid \R_+)\).

  Similarly, by \citep[Proposition~\(7.4\)d]{rtrw}, because \(\{\phi_{k,i}^<\}\) is nondecreasing with \(k\), its epi-limit is well defined and \(\phi_{k,i}^< \eto \sup_k \cl \phi_{k,i}^< = \chi(\cdot \mid \R_+)\).

  Because \(\phi_{k,i}^< \leq \phi_{k,i} \leq \phi_{k,i}^>\), \citep[Proposition~\(7.4\)g]{rtrw} implies that \(\{\phi_{k,i}\} \eto \chi(\cdot \mid \R_+)\), and consequently, we obtain \(\{\phi_k\} \eto \chi(\cdot \mid \R^n_+)\).
\end{proof}

\begin{shadytheorem}%
  \label{thm:elim-bar}
  \(\elim f + h + \phi_k = f + h + \chi(\cdot \mid \R^n_+)\).
\end{shadytheorem}

\begin{proof}
  \Cref{lem:elim-bar} and \citep[Theorem~\(7.46\)a]{rtrw} imply \(\{h + \phi_k\} \eto h + \chi(\cdot \mid \R^n_+)\).
  Finally, because \(f\) is continuous, \citep[Exercise~\(7.8\)a]{rtrw} yields \(\{f + h + \phi_k\} \eto f + h + \chi(\cdot \mid \R^n_+)\).
\end{proof}

The following corollary legitimizes the barrier approach for~\eqref{eq:nlp}.

\begin{shadycorollary}%
  \label{cor:eps-optimal}
  Let \(\inf \{f(x) + h(x) \mid x \geq 0\}\) be finite.
  For all \(\epsilon \geq 0\),
  \[
    \limsup_k (\epsilon\text{-}\argmin{}  f + h + \phi_k) \subset \epsilon\text{-}\argmin{}  f + h + \chi(\cdot \mid \R^n_+).
  \]
  In particular, if \(\{\epsilon_k\} \searrow 0\),
  \[
    \limsup_k (\epsilon_k\text{-}\argmin{} f + h + \phi_k) \subset \argmin{} f + h + \chi(\cdot \mid \R^n_+).
  \]
\end{shadycorollary}

\begin{proof}
  Follows directly from \citep[Theorem~\(7.31\)b]{rtrw}.
\end{proof}

If \(\inf \{f(x) + h(x) \mid x \geq 0\}\) is finite, the definition of the limit superior of a sequence of sets and the second part of \Cref{cor:eps-optimal} indicate that for any \(\{\epsilon_k\} \searrow 0\), there exists \(N \in \mathcal{N}_{\infty}^{\#}\) and \(\bar x_k \in \epsilon_k\text{-}\argmin{} f + h + \phi_k\) for all \(k \in N\) such that \(\{\bar x_k\}_N\) converges to a solution of~\eqref{eq:nlp}.

\subsection{Barrier subproblem}

The \(k\)-th subproblem is
\begin{equation}
  \label{eq:bar-subproblem}
  \minimize{x} \ f(x) + h(x) + \phi_k(x).
\end{equation}
\(x_k^* \in \R^n_{++}\) is first-order stationary for~\eqref{eq:bar-subproblem} if
\begin{equation}
  \label{eq:bar-subproblem-stationarity}
  0 \in \nabla f(x_k^*) - \mu_k (X_k^*)^{-1} e + \partial h(x_k^*).
\end{equation}
We call the process of solving~\eqref{eq:bar-subproblem} the \(k\)-th sequence of \emph{inner} iterations, and we denote its iterates \(x_{k,j}\) for \(j \geq 0\).
The definition of~\eqref{eq:bar-subproblem} along with certain parameter updates will be called an \emph{outer} iteration.

For \(x \in \R^n_{++}\) and \(\delta \in (0, 1)\), let
\begin{equation}
  \label{eq:def-Rn-delta-x}
  \R^n_{\delta}(x) := \{ s \in \R^n \mid \min_i (x + s)_i \geq \delta \min_i x_i \} \subset (-x + \R^n_{++}).
\end{equation}
Note that \(\R^n_{\delta}(x)\) is closed, and also convex, as shown in the following lemma.

\begin{shadylemma}
  \label{lem:rn-delta-x-convex}
  Let \(\delta \in (0,1)\) and \(x \in \R_{++}^n\).
  Then \(\R_{\delta}^n(x)\) defined in~\eqref{eq:def-Rn-delta-x} is convex.
\end{shadylemma}

\begin{proof}
  Let \(s_1\) and \(s_2 \in \R_{\delta}^n(x)\), and \(t \in [0, 1]\).
  By definition, \(\min_i (s_1 + x)_i \ge \delta \min_i x_i\) and \(\min_i (s_2 + x)_i \ge \delta \min_i x_i\).
  Now,
  \begin{align*}
    \min_i (t s_1 + (1 - t) s_2 + x)_i &\ge \min_i (t (s_1 + x))_i + \min_i ((1 - t) (s_2 + x))_i \\
    &= t \min_i (s_1 + x)_i + (1 - t)\min_i (s_2 + x)_i \\
    &\ge t \delta \min_i x_i + (1 - t) \delta \min_i x_i \\
    &= \delta \min_i x_i.
  \end{align*}
  Thus, \(t s_1 + (1 - t) s_2 \in \R_{\delta}^n(x)\) and \(\R_{\delta}^n(x)\) is convex.
\end{proof}
Under the assumptions of \Cref{lem:rn-delta-x-convex} and for \(\Delta > 0\), \(\Delta \B \cap \R_{\delta}^n(x)\) is convex.

At outer iteration \(k\), we choose \(\delta_k \in (0, 1)\), and solve~\eqref{eq:bar-subproblem} inexactly by approximately solving a sequence of trust-region subproblems of the form
\begin{subequations}
  \label{eq:tr-sub}
  \begin{align}
    \minimize{s} & \ m(s; x_{k,j}) + \chi(s \mid \Delta_{k,j} \B \cap \R^n_{\delta_k}(x_{k,j})),
    \\
    \label{eq:def-model}
    m(s; x_{k,j}) & := \varphi(s; x_{k,j}) + \psi(s; x_{k,j}),
  \end{align}
\end{subequations}
where \(\varphi(s; x_{k,j}) \approx (f + \phi_k)(x_{k,j} + s)\) and \(\psi(s; x_{k,j}) \approx h(x_{k,j} + s)\) model the smooth and nonsmooth parts of~\eqref{eq:bar-subproblem}, respectively, and \(\Delta_{k,j} > 0\) is a trust-region radius.
Models are required to satisfy the following assumption.
\begin{modelassumption}%
  \label{asm:model}
  For any \(k \in \N\) and \(j \in \N\), \(\varphi(\cdot; x_{k,j})\) is continuously differentiable on \(\R^n_{++}\) with \(\varphi(0; x_{k,j}) = f(x_{k,j}) + \phi_k(x_{k,j})\) and \(\nabla_s \varphi(0; x_{k,j}) = \nabla f(x_{k,j}) + \nabla \phi_k(x_{k,j})\).
  In addition, \(\nabla_s \varphi(\cdot; x_{k,j})\) is Lipschitz continuous with constant \(L_{k,j} \geq 0\).
  We require that \(\psi(\cdot; x_{k,j})\) be proper, lsc, and satisfy \(\psi(0; x_{k,j}) = h(x_{k,j})\) and \(\partial \psi(0; x_{k,j}) = \partial h(x_{k,j})\).
\end{modelassumption}

\begin{shadyproposition}%
  \label{prop:model-stationarity}
  Let \Cref{asm:model} be satisfied. Then \(s = 0\) is first-order stationary for~\eqref{eq:tr-sub} if and only if \(x_{k,j}\) is first-order stationary for~\eqref{eq:bar-subproblem}.
\end{shadyproposition}

\begin{proof}
  If \(x_{k,j}\) is first-order stationary, then \(0 \in \nabla f(x_{k,j}) - \mu_k X_{k,j}^{-1}e + \partial h(x_{k,j}) = \nabla_s \varphi(0; x_{k,j}) + \partial \psi(0; x_{k,j})\).
  Note that \(x_{k,j} > 0\) so \(\partial \chi(0 \mid \Delta_{k,j} \B \cap \R^n_{\delta_k}(x_{k, j})) = \{0\}\), and \(\Delta_{k,j} > 0\), so there is an open set \(O \subset \Delta_{k,j}\B \cap \R_{\delta_k}^n(x_{k,j})\) such that for all \(s \in O\), \(\chi(s \mid \Delta_{k,j} \B \cap \R^n_{\delta_k}(x_{k, j})) = 0\).
  Thus, using the definition of the subdifferential, \(\partial (\psi(\cdot; x_{k,j}) + \chi(\cdot \mid \Delta_{k,j} \B \cap \R^n_{\delta_k}(x_{k, j})))(0) = \partial \psi(0; x_{k,j})\).
  We conclude that \(0 \in \nabla_s \varphi(0; x_{k,j}) + \partial (\psi(\cdot; x_{k,j}) + \chi(\cdot \mid \Delta_{k,j} \B \cap \R^n_{\delta_k}(x_{k, j})))(0)\), which is the definition of \(s = 0\) being first-order stationary for~\eqref{eq:tr-sub}.
  The reciprocal can also be established from these observations, because if \(0 \in \nabla_s \varphi(0; x_{k,j}) + \partial (\psi(\cdot; x_{k,j}) + \chi(\cdot \mid \Delta_{k,j} \B \cap \R^n_{\delta_k}(x_{k, j})))(0)\), we have shown that \(0 \in \nabla_s \varphi(0; x_{k,j}) + \partial \psi(0; x_{k,j}) = \nabla f(x_{k,j}) - \mu_k X_{k,j}^{-1}e + \partial h(x_{k,j})\).
\end{proof}

As a special case of \Cref{prop:model-stationarity}, if \(s = 0\) solves~\eqref{eq:tr-sub}, then \(x_{k,j}\) is first-order stationary for~\eqref{eq:bar-subproblem}.

Let
\begin{equation}
  \label{eq:varphi-f}
  \varphi_f(s; x_{k, j}, B_{k,j}) := f(x_{k, j}) + \nabla f(x_{k, j})^T s + \tfrac{1}{2}s^T B_{k, j} s,
\end{equation}
where \(B_{k,j} = B_{k,j}^T\), be a second order Taylor approximation of \(f\) about \(x_{k,j}\).
We are particularly interested in the quadratic model
\begin{equation}
  \label{eq:pd-model}
  \begin{aligned}
  \varphi(s; x_{k,j}, B_{k, j}) &:=
  \varphi_f(s; x_{k, j}, B_{k,j}) + \phi_k(x_{k,j}) - \mu_k e^T X_{k,j}^{-1} s + \tfrac{1}{2} s^T X_{k,j}^{-1} Z_{k,j}s \\
  &= (f + \phi_k)(x_{k,j}) + (\nabla f(x_{k,j}) - \mu_k X_{k,j}^{-1} e)^T s + \tfrac{1}{2} s^T (B_{k,j} + X_{k,j}^{-1} Z_{k,j}) s,
  \end{aligned}
\end{equation}
where \(z_{k,j}\) is an approximation to the vector of multipliers for the bound constraints of~\eqref{eq:nlp}.

Let \(s_{k,j}\) be an approximate solution of~\eqref{eq:tr-sub}.
If \(s_{k,j}\) is accepted as a step for our algorithm used to solve~\eqref{eq:bar-subproblem} (the acceptance condition is detailed in \Cref{alg:bar-inner}), we perform the update \(x_{k,j+1} = x_{k,j} + s_{k,j}\).

By analogy with the smooth case, we use \(z_{k,j} := \mu_k X_{k,j}^{-1} e\) when \(x_{k,j}\) is first-order stationary for~\eqref{eq:bar-subproblem}.
Multiplying through by \(X_{k,j}\), we obtain \(X_{k,j} z_{k,j} = \mu_k e\).
Linearizing the continuous equality \(X z - \mu e  = 0\) with respect to \(x\) and \(z\) and evaluating all quantities at iteration \((k,j)\) yields \(X_{k,j} \Delta z_{k,j} + Z_{k,j} s_{k,j} = \mu_k e - X_{k,j} z_{k,j}\), which suggests that if \(x_{k,j+1} = x_{k,j} + s_{k,j}\), then \(z_{k,j+1} = z_{k,j} + \Delta z_{k,j} = \mu_k X_{k,j}^{-1} e - X_{k,j}^{-1} Z_{k,j} s_{k,j}\).

However, the latter \(z_{k,j+1}\) may not be positive.
We perform the update described by \citet{conn-gould-toint-2000}, by defining
\begin{equation}
  \label{eq:z-bar-update}
  \hat z_{k,j+1} = \mu_k X_{k,j}^{-1} e - X_{k,j}^{-1} Z_{k,j} s_{k,j},
\end{equation}
and projecting \(\hat z_{k,j+1}\) componentwise into the following interval to get \(z_{k,j+1}\)
\begin{equation}
  \label{eq:z-interval}
  \mathcal{I} = [\kappa_{\textup{zul}} \min (e, z_{k,j}, \mu_k X_{k, j+1}^{-1}e), \ \max (\kappa_{\textup{zuu}}e, z_{k,j}, \kappa_{\textup{zuu}} \mu_k^{-1}e, \kappa_{\textup{zuu}}\mu_k X_{k,j+1}^{-1}e)],
\end{equation}
with \(0 < \kappa_{\textup{zul}} < 1 < \kappa_{\textup{zuu}}\).
Projecting \(\hat z_{k,j+1}\) into~\eqref{eq:z-interval} always generates a positive \(z_{k,j+1}\).
The choice \(z_{k,j+1} = z_{k,j}\) is also available.
The other bounds of~\eqref{eq:z-interval} will be useful in \Cref{sec:convergence-inner} and \Cref{sec:convergence-outer}.


We define the following model, based upon a first-order Taylor approximation
\begin{subequations}%
  \begin{align}
    \varphi_{\textup{cp}}(s; x_{k,j}) &:=
    (f + \phi_k)(x_{k,j}) + (\nabla f(x_{k,j}) - \mu_k X_{k,j}^{-1}e)^T s, \label{eq:varphi-pg} \\
    m_{\textup{cp}}(s; x_{k,j}, \nu_{k,j}) &:= \varphi_{\textup{cp}}(s; x_{k,j}) + \tfrac{1}{2}\nu_{k,j}^{-1}\|s\|^2 + \psi(s; x_{k,j}), \label{eq:sk1-sub}
  \end{align}
\end{subequations}
where ``cp'' stands for ``Cauchy point''.
Let \(s_{k,j,1}\) be the solution of~\eqref{eq:tr-sub} with model \(m_{\textup{cp}}(s; x_{k,j}, \nu_{k,j})\).
As stated in \citep[Section~\(3.2\)]{aravkin-baraldi-orban-2022}, \(s_{k,j,1}\) is actually the first step of the proximal gradient method~\eqref{eq:pg-iter-general} from \(s_{k,j,0} = 0\) applied to the minimization of \(\varphi_{\textup{cp}} + \psi\) with step length \(\nu_{k,j}\):
\begin{equation}
  \label{eq:skj-1step-pg}
  s_{k,j,1} \in \prox{\nu_{k,j} \psi(\cdot; x_{k,j}) + \chi(\cdot \mid \Delta_{k,j} \B \cap \R_{\delta_k}^n(x_{k,j}))}(-\nu_{k,j}\nabla \varphi_{\textup{cp}}(0; x_{k,j}, \nu_{k,j})).
\end{equation} 
Let
\begin{equation}%
  \label{eq:optim-measure}
  \xi_{\textup{cp}}(\Delta; x_{k,j}, \nu_{k,j}) := (f + \phi_k + h)(x_{k,j}) - (\varphi_{\textup{cp}} + \psi)(s_{k,j,1}; x_{k,j}),
\end{equation}
where \(\Delta > 0\), and let \(\nu_{k,j}^{-1/2} \xi_{\textup{cp}}(\Delta; x_{k,j}, \nu_{k,j})^{1/2}\) be our measure of criticality.
\citet{aravkin-baraldi-orban-2022} and its corrigendum \citep{aravkin-baraldi-leconte-orban-2023} indicate that \(\nu_{k,j}^{-1/2} \xi_{\textup{cp}}(\Delta; x_{k,j}, \nu_{k,j})^{1/2}\) is similar to \(\nu_{k,j}^{-1} \|s_{k,j,1}\|\), which is the norm of the generalized gradient at \(x_{k,j}\).
We can apply \citep[Theorems~\(1.17\) and~\(7.41\)]{rtrw} to conclude that \(\xi_{\textup{cp}}(\Delta; x_{k,j}, \nu_{k,j})\) is proper lsc in \((x_{k,j}, \nu_{k,j}) \in \R_{++}^n \times \R_{++}\).
In particular, \(\nu_{k,j}^{-1/2} \xi_{\textup{cp}}(\Delta; x_{k,j}, \nu_{k,j})^{1/2} = 0\) for any \(\Delta > 0\) and \(\nu_{k,j} > 0\) \(\Longrightarrow\) \(s = 0\) solves~\eqref{eq:tr-sub}, and \(x_{k,j}\) is first-order stationary for~\eqref{eq:bar-subproblem}.

\Cref{alg:bar-outer} summarizes the outer iteration.

\begin{algorithm}[h]
  \caption[caption]{%
    Nonsmooth interior-point method (outer iteration).%
    \label{alg:bar-outer}
  }
  \begin{algorithmic}[1]
    \State Choose \(\epsilon > 0\), sequences \(\{\mu_k\} \searrow 0\), \(\{\epsilon_{d,k}\} \searrow 0\), \(\{\epsilon_{p,k}\} \searrow 0\), and \(\{\delta_k\} \rightarrow \bar \delta \in [0,1)\) with \(\delta_k \in (0,1)\) for all \(k\).
    \State Choose \(x_{0, 0} \in \R^n_{++}\) where \(h\) is finite.
    \For{\(k = 0, 1, \ldots\)}
      \State Compute an approximate solution \(x_k := x_{k,j}\) to~\eqref{eq:bar-subproblem} and \(z_k := z_{k,j}\) in the sense that
      \begin{equation}
        \label{eq:stat-inner-subproblem}
        \nu_{k,j}^{-1/2} \xi_{\textup{cp}}(\Delta_{k,j}; x_{k,j}, \nu_{k,j})^{1/2} \leq \epsilon_{d,k}
      \end{equation}
      and 
      \begin{equation}
        \label{eq:central-path-inner}
        \|X_{k,j} z_{k,j} - \mu_k e\| \le \epsilon_{p,k}.
      \end{equation}
      \State Set \(x_{k+1, 0} := x_k\).
    \EndFor
  \end{algorithmic}
\end{algorithm}

For each outer iteration \(k\), the inner iterations generate a sequence \(\{x_{k,j}\}\) according to an adaptation of \citep[Algorithm~\(3.1\)]{aravkin-baraldi-orban-2022} in which the subproblems have the form~\eqref{eq:tr-sub} with the smooth part of the model defined by~\eqref{eq:pd-model}.
Each trust-region step is required to satisfy the following assumption.

\begin{stepassumption}%
  \label{asm:cauchy-decrease}
  Let \(k \in \N\).
  There exists \(\kappa_{\textup{m},k} > 0\) and \(\kappa_{\textup{mdc},k} \in (0, \, 1)\) such that for all \(j\), \(s_{k,j} \in \Delta_{k,j} \B \cap \R_{\delta_k}(x_{k,j})\),
  \begin{subequations}
    \begin{align}
    \label{eq:model-adequation}
      |(f + \phi_k + h)(x_{k,j} + s_{k,j}) - m(s_{k,j}; x_{k,j}, B_{k,j})| & \leq \kappa_{\textup{m},k} \|s_{k,j}\|^2,
      \\
      \label{eq:cauchy-decrease}
      m(0; x_{k,j}, B_{k,j}) - m(s_{k,j}; x_{k,j}, B_{k,j}) & \geq \kappa_{\textup{mdc},k} \xi_{\textup{cp}}(\Delta_{k,j}; x_{k,j}, \nu_{k,j}),
    \end{align}
  \end{subequations}
  where \(m\) is defined in~\eqref{eq:def-model}--\eqref{eq:pd-model}, and \(\xi_{\textup{cp}}(\Delta_{k,j}; x_{k,j}, \nu_{k,j})\) is defined in~\eqref{eq:optim-measure}.
\end{stepassumption}

In \Cref{asm:cauchy-decrease}, the subscript ``m'' of \(\kappa_{\textup{m},k}\) refers to the model adequacy, and the subscript ``mdc'' of \(\kappa_{\textup{mdc},k}\) refers to the model decrease. 

\Cref{alg:bar-inner} summarizes the process.

\begin{algorithm}[h]
  \caption[caption]{%
    Nonsmooth interior-point method (inner iteration).%
    \label{alg:bar-inner}
  }
  \begin{algorithmic}[1]
    \State Choose constants
    \[
      0 < \eta_1 \leq \eta_2 < 1,
      \mathhfill
      0 < \gamma_1 \leq \gamma_2 < 1 < \gamma_3 \leq \gamma_4,
      \mathhfill
      \Delta_{\max} > \Delta_{k,0} > 0,
      \mathhfill
      \alpha > 0, \, \beta \geq 1.
    \]
    \State Compute \(f(x_{k,0}) + \phi_k(x_{k,0}) + h(x_{k,0})\).
    \For{\(j = 0, 1, \ldots\)}
      \State\label{alg:tr-nonsmooth:nuk}%
      Choose \(0 < \nu_{k,j} \leq  1 / (L_{k,j} + \alpha^{-1} \Delta_{k,j}^{-1})\).
      \State Define \(m(s; x_{k,j}, B_{k,j})\) as in~\eqref{eq:def-model} satisfying \Cref{asm:model}.
      \State Define \(m_{\textup{cp}}(s; x_{k,j}, \nu_{k,j})\) as in~\eqref{eq:sk1-sub}.
      \State\label{alg:tr-nonsmooth:sk1}%
      Compute \(s_{k,j,1}\) as a solution of~\eqref{eq:skj-1step-pg}.
      \State\label{alg:tr-nonsmooth:sk}%
      Compute an approximate solution \(s_{k,j}\) of~\eqref{eq:tr-sub} such that \(\|s_{k,j}\| \leq \min(\Delta_{k,j}, \, \beta \|s_{k,j,1}\|)\).
      \State Compute the ratio
      \[
      \rho_{k,j} :=
      \frac{
        (f + \phi_k + h)(x_{k,j}) - (f + \phi_k + h)(x_{k,j} + s_{k,j})
      }{
        m(0; x_{k,j}, B_{k,j}) - m(s_{k,j}; x_{k,j}, B_{k,j})
      }.
      \]
      \State If \(\rho_{k,j} \geq \eta_1\), set \(x_{k,j+1} = x_{k,j} + s_{k,j}\) and update \(z_{k,j+1}\) according to~\eqref{eq:z-bar-update} and~\eqref{eq:z-interval}.
      Otherwise, set \(x_{k,j+1} = x_{k,j}\) and \(z_{k,j+1} = z_{k,j}\).
      \State Update the trust-region radius according to
      \[
        \hat \Delta_{k,j+1} \in
        \left\{
          \begin{array}{lll}
             {[\gamma_3 \Delta_{k,j}, \, \gamma_4 \Delta_{k,j}]} &
             \text{ if } \rho_{k,j} \geq \eta_2, &
             \text{(very successful iteration)}
          \\ {[\gamma_2 \Delta_{k,j}, \, \Delta_{k,j}]} &
             \text{ if } \eta_1 \leq \rho_{k,j} < \eta_2, &
             \text{(successful iteration)}
          \\ {[\gamma_1 \Delta_{k,j}, \, \gamma_2 \Delta_{k,j}]} &
             \text{ if } \rho_{k,j} < \eta_1 &
             \text{(unsuccessful iteration)}
          \end{array}
        \right.
      \]
      and \(\Delta_{k,j+1} = \min (\hat \Delta_{k,j+1}, \, \Delta_{\max})\).
    \EndFor
  \end{algorithmic}
\end{algorithm}

\subsection{Convergence of the inner iterations}
\label{sec:convergence-inner}


Let \(k\) and \(j\) be fixed positive integers, and \(x_{k, j} > 0\).
We may rewrite our ``Cauchy point'' subproblem for the inner iterations as in \citep{aravkin-baraldi-orban-2022}:
\begin{subequations}
  \begin{align}
    \label{eq:def-model-sub}
    p(\Delta; x_{k, j}, \nu_{k,j}, \delta_k) &:= \minimize{s} m_{\textup{cp}}(s; x_{k,j}, \nu_{k,j}) + \chi (s; \Delta \B \cap \R_{\delta_k}^n(x_{k, j})) \\
    P(\Delta; x_{k, j}, \nu_{k,j}, \delta_k) &:= \argmin{s} m_{\textup{cp}}(s; x_{k,j}, \nu_{k,j}) + \chi (s; \Delta \B \cap \R_{\delta_k}^n(x_{k, j})).
  \end{align}
\end{subequations}

First, we present some properties of the subproblem~\eqref{eq:def-model-sub} in the following result.

\begin{shadyproposition}[{\protect \citealp[Proposition~\(3.1\)]{aravkin-baraldi-orban-2022}}]
  \label{prop:sub-model}
  Let \Cref{asm:model} be satisfied, \(\nu > 0\), \(\delta > 0\) and \(x \in \R_+^n\).
  If we define \(p(0; x, \nu, \delta) := \varphi_{\textup{cp}}(0; x) + \psi(0; x)\) and \(P(0; x, \nu, \delta) = \{0\}\), the domain of \(p(\cdot; x, \nu, \delta)\) and \(P(\cdot; x, \nu, \delta)\) is \(\{\Delta \mid \Delta \ge 0\}\).
  In addition,
  \begin{enumerate}
    \item \(p(\cdot; x, \nu, \delta)\) is proper lsc and for each \(\Delta \ge 0\), \(P(\Delta; x, \nu, \delta)\) is nonempty and compact;
    \item if \(\{\Delta_{k,j}\} \rightarrow \bar \Delta_k \ge 0\) in such a way that \(\{p(\Delta_{k, j}; x, \nu, \delta)\} \rightarrow p(\bar \Delta_k; x, \nu, \delta)\), and for each \(j\), \(s_{k, j} \in P(\Delta_{k, j}; x, \nu, \delta)\), then \(\{s_{k, j}\}\) is bounded and all its limit points are in \(P(\bar \Delta_k; x, \nu, \delta)\);
    \item if \(\varphi_{\textup{cp}}(\cdot; x) + \tfrac{1}{2}\nu^{-1}\|s\|^2 + \psi(\cdot; x)\) is strictly convex, \(P(\Delta; x, \nu, \delta)\) is single valued;
    \item if \(\bar \Delta_k > 0 \) and there exists \(\bar s \in P(\bar \Delta_k; x, \nu, \delta)\) such that \(\bar s \in \mathrm{int}( \bar \Delta_k \B \cap \R_{\delta_k}^n(x))\), then \(p(\cdot; x, \nu, \delta)\) is continuous at \(\bar \Delta_k\) and \(\{p(\Delta_{k, j}; x, \nu, \delta)\} \rightarrow p(\bar \Delta_k; x, \nu, \delta)\) holds in part 2. 
  \end{enumerate}
\end{shadyproposition}
\begin{proof}
  \Cref{asm:model} and the compactness of \(\Delta \B \cap \R_{\delta_k}^n (x)\) ensure that the objective of~\eqref{eq:def-model-sub} is always level-bounded in \(s\) locally uniformly in \(\Delta\), because for any \(\bar \Delta > 0\), \(\epsilon > 0\), and \(\Delta \in (\bar \Delta - \epsilon, \bar \Delta + \epsilon)\) with \(\Delta \ge 0\), its level sets are contained in \(\Delta \B \cap \R_{\delta_k}(x_{k, j}) \subseteq (\bar \Delta + \epsilon) \B \cap \R_{\delta_k}(x_{k, j})\).
  From this observation, we can draw similar conclusions to the analysis of \citep[Proposition~\(3.1\)]{aravkin-baraldi-orban-2022}.
\end{proof}

The observation in the proof of \Cref{prop:sub-model} and \Cref{asm:model} allows us to derive directly some of the convergence properties of \citet{aravkin-baraldi-orban-2022} for \Cref{alg:bar-inner}.

\begin{shadyproposition}[{\protect \citealp[Theorem~\(3.4\)]{aravkin-baraldi-orban-2022}}]
  \label{prop:delta-succ}
  Let \Cref{asm:model} and \Cref{asm:cauchy-decrease} be satisfied and let
  \begin{equation}
    \Delta_{\textup{succ},k} := \frac{\kappa_{\textup{mdc},k} (1 - \eta_2)}{2 \kappa_{\textup{m},k} \alpha \beta^2} > 0.
  \end{equation}
  If \(x_{k, j}\) is not first-order stationary for~\eqref{eq:bar-subproblem} and \(\Delta_{k, j} \le \Delta_{\textup{succ}}\), then iteration \(j\) is very successful and \(\Delta_{k, j+1} \ge \Delta_{k, j}\).
\end{shadyproposition}
\begin{proof}
  If \(x_{k,j}\) is not first-order stationary, \(\R_{\delta_k}^n(x_{k, j}) \neq \{0\}\), thus \(s_{k,j,1} \neq 0\) and \(s_{k,j} \neq 0\).
  The rest of the proof is identical to that of \citep[Theorem~\(3.4\)]{aravkin-baraldi-orban-2022}.
\end{proof}

Now, let \(\Delta_{\min, k} := \min (\Delta_{k,0}, \gamma_1 \Delta_{\textup{succ},k}) > 0\).
Then, \(\Delta_{k,j} \ge \Delta_{\min,k}\) for all \(j \in \N\).
If we consider \(\varphi\) defined in~\eqref{eq:pd-model}, for \(s_1\) and \(s_2\) in \(\R_+^n\),
\begin{equation}
  \label{eq:varphi-grad-diff}
  \|\nabla \varphi(s_1) - \nabla \varphi(s_2) \| = \|B_{k, j}(s_1 - s_2) + X_{k, j}^{-1} Z_{k, j} (s_1 - s_2)\|.
\end{equation}
As \(L_{k,j} = \|B_{k,j} + X_{k,j}^{-1}Z_{k,j}\| \le \|B_{k,j}\| + \|X_{k,j}^{-1}\|\|Z_{k,j}\|\), if \(\{B_{k,j}\}_j\) remains bounded, \(\{x_{k,j}\}_j\) must be bounded away from zero to guarantee the existence of some \(M_k > 0\) such that \(L_{k,j} \le M_k\) for all \(j \in \N\).
To apply the complexity results, and to establish that \(\liminf \nu_{k, j}^{-1/2} \xi_{\textup{cp}}(\Delta_{k,j}; x_{k, j}, \nu_{k,j})^{1/2} = 0\) if \(f + h + \phi_k\) is bounded below on \(\R_+^n\), we need a stronger assumption on the Lipschitz constant of our model.

\begin{modelassumption}[{\protect \citealp[Model Assumption~\(3.3\)]{aravkin-baraldi-orban-2022}}]
  \label{asm:lipschitz-bnd}
  In \Cref{asm:model}, there exists \(M_k > 0\) such that \(0 \le L_{k, j} \le M_k\) for all \(j \in \N\).
  In addition, we select \(\nu_{k, j}\) in line 4 of \Cref{alg:bar-inner} in a way that there exists \(\nu_{\min, k} > 0\) such that \(\nu_{k, j} \ge \nu_{\min, k}\) for all \(j \in \N\).
\end{modelassumption}

As in \citet{aravkin-baraldi-orban-2022}, we can set \(\nu_{k, j} := 1 / (L_{k, j} + \alpha^{-1}\Delta_{k, j}^{-1})\) in \Cref{alg:bar-inner} to ensure \(\nu_{k, j} \ge \nu_{\min, k} := 1 / (M_k + \alpha^{-1} \Delta_{\min, k}^{-1}) > 0\) if the first part of \Cref{asm:lipschitz-bnd} holds with \(M_k\).
The observations below~\eqref{eq:varphi-grad-diff} motivate us to prove that \(\{x_{k,j}\}_j\) is bounded away from zero in the next result.

\begin{shadyproposition}
  \label{prop:min-xkji-kappa-mdb}
  Let \(k \in \N\), \Cref{asm:model} be satisfied for \(\varphi\) in~\eqref{eq:pd-model}, and \((f + h)(x_{k, j}) \ge (f + h)_{\textup{low}, k}\).
  Then, there exists \(\kappa_{\textup{mdb}, k} > 0\) such that, for all \(j\), we have
  \begin{equation}
    \label{eq:const-dist-xkj}
    \min_i (x_{k, j})_i \ge \kappa_{\textup{mdb}, k}.
  \end{equation}
\end{shadyproposition}
\begin{proof}
  We proceed similarly as in \citet[Theorem~\(13.2.1\)]{conn-gould-toint-2000}.
  Let \(k\) be a positive integer and \(\{x_{k, j}\}\) be a sequence generated by \Cref{alg:bar-inner}.
  As \(\{(f + \phi_k + h)(x_{k,j})\}_j\) is decreasing and \((f + h)(x_{k, j}) \ge (f + h)_{\textup{low}, k}\), we have \(\limsup_j \phi_k(x_{k,j}) < \infty\), which implies that~\eqref{eq:const-dist-xkj} holds.
\end{proof}

The following proposition shows that we can use the convergence results of \citet{aravkin-baraldi-orban-2022} using the same assumptions they used for \(\varphi_f\).
It will justify that \Cref{asm:lipschitz-bnd} can be used for \(\varphi\) defined in~\eqref{eq:pd-model}.

\begin{shadyproposition}
  \label{prop:varphi-lipschitz-grad}
  Under the assumptions of \Cref{prop:min-xkji-kappa-mdb}, let \(\varphi_f\) be defined as in~\eqref{eq:varphi-f} so that \(\nabla_s \varphi_f(\cdot; x_{k,j}, B_{k,j})\) is Lipschitz continuous with constant \(\tilde L_{k,j} \geq 0\) and there exists \(\tilde M_k > 0\) such that \(0 \le \tilde L_{k,j} \le \tilde M_k\) for all \(j \in \N\).
  Then \(\varphi\) satisfies \Cref{asm:lipschitz-bnd}.
\end{shadyproposition}
\begin{proof}
  We can use~\eqref{eq:z-interval} and~\eqref{eq:const-dist-xkj} to say that \(X_{k, j}^{-1} Z_{k, j}\) is bounded for all \(j\), and we deduce from~\eqref{eq:varphi-grad-diff} that \(\varphi\) satisfies \cref{asm:lipschitz-bnd}.
\end{proof}

Now, we justify that \Cref{asm:cauchy-decrease} holds when \(h(x_{k,j} + s_{k,j}) = \psi(s_{k,j}; x_{k,j})\).
As \(\{x_{k,j}\}_j\) remains bounded away from \(\partial \R_+^n\) with \Cref{prop:min-xkji-kappa-mdb}, so does \(\{x_{k,j} + s_{k,j}\}_j\) by definition of \(\R_{\delta_k}(x_{k,j})\).
Since \(\nabla f\) is Lipschitz-continuous,
\begin{equation}
  \label{eq:f-descent-ineq}
  f(x_{k,j} + s_{k,j}) - f(x_{k,j}) - \nabla f(x_{k,j})^T s_{k,j} \le \tfrac{1}{2} L_f \|s_{k,j}\|^2,
\end{equation}
and a second-order Taylor approximation of \(\phi_k\) about \(x_{k,j} + s_{k,j}\) gives 
\begin{equation}
  \label{eq:phik-taylor}
  \phi_k(x_{k,j} + s_{k,j}) - \phi_k(x_{k,j}) - \mu_k s_{k,j}^T X_{k,j}^{-1}e = 
  \tfrac{1}{2} \mu_k s_{k,j}^T X_{k,j}^{-2}s_{k,j} + o(\|s_{k,j}\|^2). 
\end{equation}
Under the assumptions of \Cref{prop:varphi-lipschitz-grad}, \(L_{k,j} = \|B_{k,j} + X_{k,j}^{-1}Z_{k,j}\| \le M_k\), and
\begin{multline*}
  |(f + \phi_k + h)(x_{k,j} + s_{k,j}) - m(s_{k,j}; x_{k,j}, B_{k,j})| \\
    = |f(x_{k,j} + s_{k,j}) - f(x_{k,j}) - \nabla f(x_{k,j})^T s_{k,j} + \phi_k(x_{k,j} + s_{k,j}) - \phi_k(x_{k,j}) - \\
    \mu_k s_{k,j}^T X_{k,j}^{-1}e - \tfrac{1}{2}s_{k,j}^T(B_{k,j} + X_{k,j}^{-1}Z_{k,j})s_{k,j}| + o(\|s_{k,j}\|^2).
\end{multline*}
The above equality combined with~\eqref{eq:f-descent-ineq} and~\eqref{eq:phik-taylor} implies that~\eqref{eq:model-adequation} holds.

To emphasize the similarities between our inner iterations and the trust-region algorithm of \citet{aravkin-baraldi-orban-2022}, and in light of \Cref{prop:varphi-lipschitz-grad}, we use in our next results that \(\varphi\) satisfies \Cref{asm:lipschitz-bnd}, instead of writting assumptions on \(\varphi_f\).
The following proposition gives us a sufficient condition for~\eqref{eq:cauchy-decrease} to be satisfied.

\begin{shadyproposition}[{\protect \citealp[Proposition~\(1\)]{aravkin-baraldi-leconte-orban-2023}}]
  \label{prop:justif-model-xi-ineq}
  If \Cref{asm:lipschitz-bnd} is satisfied with bounded Hessian approximations \(\{B_{k,j}\}_j\), then there exists \(\kappa_{\textup{mdc},k} \in (0,1)\) such that~\eqref{eq:cauchy-decrease} holds for all \(j\).
\end{shadyproposition}

\begin{proof}
  The proof is identical to that of \citep[Proposition~\(1\)]{aravkin-baraldi-leconte-orban-2023} when replacing \(B_k\) by \(B_{k,j} + X_{k,j}^{-1}Z_{k,j}\), and using the subscripts\(~_{k,j}\) where \(j\) is the iteration number of the algorithm instead of the subscript\(~_k\).
\end{proof}

\Cref{prop:min-xkji-kappa-mdb} allows us to write the following convergence results for \cref{alg:bar-inner}.
As in \citet{aravkin-baraldi-orban-2022}, we define the smallest iteration number \(j_k(\epsilon)\) such that 
\begin{equation}
  \label{eq:stat-inner}
  \nu_{k, j}^{-1/2}\xi_{\textup{cp}}(\Delta_{k,j}; x_{k,j}, \nu_{k,j})^{1/2} \le \epsilon \quad (0 < \epsilon < 1), \quad k = 0, 1, \ldots,
\end{equation}
and we express the set of successful iterations, the set of successful iterations for which~\eqref{eq:stat-inner} has not yet been attained, and the set of unsuccessful iterations for which~\eqref{eq:stat-inner} has not yet been attained as
\begin{subequations}
  \begin{align}
    \mathcal{S}_k &:= \{j \in \N \mid \rho_{k, j} \ge \eta_1\} \\
    \mathcal{S}_k(\epsilon) &:= \{j \in \mathcal{S}_k \mid j < j_k(\epsilon)\} \\
    \mathcal{U}_k(\epsilon) &:= \{j \in \N \mid j \notin \mathcal{S}_k \text{ and } j < j_k(\epsilon)\}.
  \end{align}
\end{subequations}

\begin{shadyproposition}[{\protect \citealp[Theorem~\(3.5\)]{aravkin-baraldi-orban-2022}}]
  Let \Cref{asm:model} and \Cref{asm:cauchy-decrease} be satisfied.
  If \Cref{alg:bar-inner} only generates finitely many successful iterations, then \(x_{k, j} = x_k^*\) for all sufficiently large \(j\) and \(x_k^*\) is first-order critical for~\eqref{eq:bar-subproblem}.
\end{shadyproposition}
\begin{proof}
  The proof is inspired from \citep[Theorem~\(3.5\)]{aravkin-baraldi-orban-2022}, which itself follows that of \citet[Theorem~\(6.6.4\)]{conn-gould-toint-2000}.
  The assumptions indicate that there is \(j_0 \in \N\) such that all iterations \(j \ge j_0\) are unsuccessful, and \(x_{k,j} = x_{k,j_0} = x_k^*\) because of the update rules of \Cref{alg:bar-inner}.
  We assume by contradiction that \(x_k^*\) is not first-order critical.
  \(x_k^*\) does not have any of its components equal to \(+\infty\) because it is attained after a finite number of iterations of \Cref{alg:bar-inner}.
  As \(h\) is proper, \((f + h)(x_k^*) > - \infty\).
  Thus, \Cref{prop:min-xkji-kappa-mdb} implies that there exists \(\Delta_k^*\) such that \(\Delta_k^* \B \subset \R_{\delta_k}(x_k^*)\).
  Since all iterations \(j \ge j_0\) are unsuccessful, there will be some \(j_1 \ge j_0\) such that \(\Delta_{j_1} \le \min(\Delta_{\textup{succ},k}, \, \Delta_k^*)\), which implies that iteration \(j_1\) is very successful with \Cref{prop:delta-succ}, and contradicts the fact that \(x_k^*\) is not first-order critical.
\end{proof}

Finally, we have the following result for the inner iterates.
\begin{shadyproposition}[{\protect \citealp[Theorem~\(3.11\)]{aravkin-baraldi-orban-2022}}]
  \label{prop:cv-inner}
  Let \(k \in \N\), \Cref{asm:model} and \Cref{asm:lipschitz-bnd} be verified for \(\varphi\) defined in~\eqref{eq:pd-model}, and \Cref{asm:cauchy-decrease} be satisfied for \cref{alg:bar-inner}.
  If there are infinitely many successful iterations, then, either
  \begin{equation}
    \underset{j \rightarrow \infty}{\lim} (f + \phi_k + h)(x_{k, j}) = -\infty \quad \text{or} \quad \underset{j \rightarrow \infty}{\lim} \nu_{k, j}^{-1/2} \xi_{\textup{cp}}(\Delta_{k,j}; x_{k,j}, \nu_{k,j})^{1/2} = 0.
  \end{equation}
\end{shadyproposition}
\begin{proof}
  The proof is identical as that of \citep[Theorem~\(3.11\)]{aravkin-baraldi-orban-2022}.
\end{proof}

Now, the definitions of~\(\varphi_{\textup{cp}}\) in~\eqref{eq:varphi-pg} and \(s_{k,j,1}\) as a minimizer of~\eqref{eq:tr-sub} with model \(m_{\textup{cp}}\) indicate that
\begin{multline*}
  (f + \phi_k + h)(x_{k,j}) = \varphi_{\textup{cp}}(0; x_{k,j}) + \psi(0; x_{k,j}) = m_{\textup{cp}}(0; x_{k,j}, \nu_{k,j}) \ge \\
   m_{\textup{cp}}(s_{k,j,1}; x_{k,j}, \nu_{k,j}) = \varphi_{\textup{cp}}(s_{k,j,1}; x_{k,j}) + \tfrac{1}{2}\nu_{k,j}^{-1}\|s_{k,j,1}\|^2 + \psi(s_{k,j,1}; x_{k,j}).
\end{multline*}
By reinjecting this inequality into the definition of~\(\xi_{\textup{cp}}\) in~\eqref{eq:optim-measure}, we obtain
\begin{equation*}
  \xi_{\textup{cp}}(\Delta_{k,j}; x_{k,j}, \nu_{k,j}) \ge \tfrac{1}{2} \nu_{k,j}^{-1} \|s_{k,j,1}\|^2,
\end{equation*}
so that
\begin{equation}
  \label{eq:sqrt-nu-xicp-ineq-skj1}
  \nu_{k,j}^{-1/2} \xi_{\textup{cp}}(\Delta_{k,j}; x_{k,j}, \nu_{k,j})^{1/2} \ge \tfrac{1}{\sqrt{2}} \nu_{k,j}^{-1} \|s_{k,j,1}\|.
\end{equation}
From this observation, we deduce that \(s_{k,j} \underset{j \rightarrow \infty}{\longrightarrow} 0\) in the following result.
\begin{shadylemma}
  \label{lem:skj-0}
  If \cref{asm:model} holds for \(\varphi\) defined in~\eqref{eq:pd-model} and
  \begin{equation*}
    \lim_{j \rightarrow \infty} \nu_{k,j}^{-1/2}\xi_{\textup{cp}}(\Delta_{k,j}; x_{k,j}, \nu_{k,j})^{1/2} \rightarrow 0,
  \end{equation*}
  then \(\lim_{j\rightarrow \infty} \|s_{k,j,1}\| = \lim_{j\rightarrow \infty} \|s_{k,j}\| = 0\). 
\end{shadylemma}
\begin{proof}
  We use \(\beta \|s_{k,j,1}\|^2 \ge \|s_{k,j}\|^2\) and~\eqref{eq:sqrt-nu-xicp-ineq-skj1} to conclude that
  \begin{equation*}
    \nu_{k,j}^{-1/2} \xi_{\textup{cp}}(\Delta_{k,j}; x_{k,j}, \nu_{k,j})^{1/2} \ge \tfrac{1}{\sqrt{2}} \nu_{k,j}^{-1} \|s_{k,j,1}\| \ge \tfrac{1}{\sqrt{2}} \nu_{k,j}^{-1} \beta^{-1/2} \|s_{k,j}\| \ge 0.
  \end{equation*}
  With \Cref{asm:model}, \(\nu_{k,j} \underset{j \rightarrow \infty}{\longrightarrow} \bar \nu_k > 0\), and we have \(\|s_{k,j}\| \underset{j \rightarrow \infty}{\longrightarrow} 0\).
\end{proof}

Now, we study the asymptotic satisfaction of the inner perturbed complementarity.
We show that~\eqref{eq:central-path-inner} is eventually satisfied, similarly to \citep[Theorem~\(13.6.4\)]{conn-gould-toint-2000} in the smooth case.

\begin{shadyproposition}
  \label{prop:complementarity}
  Let \(k \in \N\), \Cref{asm:model} and \Cref{asm:lipschitz-bnd} be verified for \(\varphi\) defined in~\eqref{eq:pd-model}, \Cref{asm:cauchy-decrease} be satisfied for \cref{alg:bar-inner}, and \((f + \phi_k + h)(x_{k, j}) \ge (f + \phi_k + h)_{\textup{low}, k}\) for all \(j \in \N\).
  Then,
  \begin{equation*}
    \lim_{j \rightarrow \infty} \|\mu_k X_{k, j}^{-1}e - z_{k, j}\| = 0.
  \end{equation*}
\end{shadyproposition}
\begin{proof}
  We proceed similarly as in the proof of \citep[Theorem~\(13.6.4\)]{conn-gould-toint-2000}.
  With the formula \(\hat z_{k,j+1} = \mu_k X_{k,j}^{-1}e - X_{k,j}^{-1}Z_{k,j}s_{k,j}\), 
  \begin{equation*}
    \begin{aligned}
    \|\hat z_{k,j+1} - \mu_k X_{k,j+1}^{-1}e\| &\le \|X_{k,j}^{-1} Z_{k,j} s_{k,j}\| + \mu_k \|X_{k,j+1}^{-1}e - X_{k,j}^{-1}e\|\\
    &\le \|X_{k,j}^{-1} Z_{k,j}\| \|s_{k,j}\| + \mu_k \sqrt{n} \|X_{k,j+1}^{-1} - X_{k,j}^{-1}\|. 
    \end{aligned}
  \end{equation*}
  Using \Cref{prop:cv-inner}, we have \(\lim_{j\rightarrow \infty} \nu_{k,j}^{-1/2} \xi_{\textup{cp}}(\Delta_{k,j}; x_{k,j}, \nu_{k,j})^{1/2} = 0\).
  This leads to \(\lim_{j \rightarrow \infty} \|s_{k,j}\| = 0\) with \Cref{lem:skj-0}, and we have that \(x_{k, j}\) is bounded away from \(0\) using \Cref{prop:min-xkji-kappa-mdb}.
  Therefore, either iteration \(j\) is not successful and \(X_{k,j+1} = X_{k,j}\), or
  \begin{align*}
    \|X_{k,j+1}^{-1} - X_{k,j}^{-1}\| & = \|X_{k,j}^{-1} (X_{k,j}X_{k,j+1}^{-1} - I)\|\\
                                      & = \|X_{k,j}^{-1} (X_{k,j}(X_{k,j} + S_{k,j})^{-1} - I)\|\\
                                      & = \|X_{k,j}^{-1} ((I + X_{k,j}^{-1}S_{k,j})^{-1} - I)\| \rightarrow 0.
  \end{align*}
  Since \(X_{k,j+1}^{-1} z_{k,j}\) is also bounded for all \(j\) using \Cref{prop:min-xkji-kappa-mdb} and~\eqref{eq:z-interval}, we have
  \begin{equation}
    \label{eq:lim-bar-z-mu-xinv}
    \lim_{j \rightarrow \infty}\|\hat z_{k,j+1} - \mu_k X_{k,j+1}^{-1}e\| \rightarrow 0.
  \end{equation}
  For \(j\) large enough, we have
  \begin{equation}
    \label{eq:bar-z-interval-xinv-mu}
    \kappa_{\textup{zul}} \mu_k X_{k,j+1}^{-1} e \le \hat z_{k,j+1} \le \kappa_{\textup{zuu}} \mu_k X_{k,j+1}^{-1}e,
  \end{equation}
  so that \(\hat z_{k,j+1} = z_{k,j+1}\) if \(j\) is large enough.
\end{proof}

In \Cref{alg:bar-outer}, each subproblem is solved approximately with tolerances \(\epsilon_{d,k} \searrow 0\) and \(\epsilon_{p,k} \searrow 0\).
This gives rise to the analysis in \Cref{sec:convergence-outer}.

\subsection{Convergence of the outer iterations}
\label{sec:convergence-outer}

For each \(k \in \N\), the stopping condition of \Cref{alg:bar-inner} occurs in a finite number of iterations.
Let \(j_k\) denote the number of iterations performed by \Cref{alg:bar-inner} at outer iteration \(k\).
To simplify the notation, let \(\bar x_k := x_{k,j_k}\), \(\bar z_k := z_{k,j_k}\), \(\bar s_{k,1} := s_{k,j_k,1}\), \(\bar \Delta_k := \Delta_{k,j_k}\) and \(\bar \nu_{k} := \nu_{k,j_k}\).

First, we present two assumptions that will be useful for our analysis.
\begin{parameterassumption}
  \label{asm:liminf-delta-nu}
  \(\liminf \Delta_{\min, k} = \bar \Delta > 0\)
\end{parameterassumption}

To satisfy the second part of \Cref{asm:liminf-delta-nu}, we need to change \(\varphi\) in~\eqref{eq:pd-model} to
\begin{equation}
  \label{eq:pd-model-theta}
  \begin{aligned}
    \varphi(s; x_{k,j}, B_{k, j}) &:=
    \varphi_f(s; x_{k, j}, B_{k,j}) + \phi_k(x_{k,j}) - \mu_k e^T X_{k,j}^{-1} s + \tfrac{1}{2} s^T \Theta_{k,j} s \\
    &= (f + \phi_k)(x_{k,j}) + (\nabla f(x_{k,j}) - \mu_k X_{k,j}^{-1} e)^T s + \tfrac{1}{2} s^T (B_{k,j} + \Theta_{k,j}) s,
  \end{aligned}
\end{equation}
where \(\Theta_{k,j} = \min(X_{k,j}^{-1} Z_{k,j}, \kappa_{\textup{bar}} I)\) with the \(\min\) taken componentwise and \(\kappa_{\textup{bar}} > 0\).

Since the results of \Cref{sec:convergence-inner} involving \(\varphi\) defined in~\eqref{eq:pd-model} are all based upon \Cref{asm:cauchy-decrease}, those results continue to apply if \Cref{asm:cauchy-decrease} holds for \(\varphi\) defined in~\eqref{eq:pd-model-theta}.
We now show that that is the case.

Fist, we observe that
\begin{multline*}
  (f + \phi_k)(x_{k,j} + s_{k,j}) - \varphi(s_{k,j}; x_{k,j}, B_{k,j}) = f(x_{k,j} + s_{k,j}) - f(x_{k,j}) - \nabla f(x_{k,j})^T s_{k,j} + \\
    \phi_k(x_{k,j} + s_{k,j}) - \phi_k(x_{k,j}) - \mu_k s_{k,j}^T X_{k,j}^{-1}e - \tfrac{1}{2}s_{k,j}^T(B_{k,j} + \Theta_{k,j})s_{k,j}.
\end{multline*}
Assume, as in \Cref{prop:varphi-lipschitz-grad}, that \(\|B_{k,j}\| \le \tilde M_k\).
Since \(\Theta_{k,j}\) is bounded by definition, we use~\eqref{eq:f-descent-ineq} and~\eqref{eq:phik-taylor} to conclude that \(f(x_{k,j} + s_{k,j}) + \phi_k(x_{k,j} + s_{k,j}) - \varphi(s_{k,j}; x_{k,j}, B_{k,j}) = O(\|s_{k,j}\|^2)\), and, if \(\psi(s; x_{k,j}) = h(x_{k,j} + s)\),~\eqref{eq:model-adequation} holds.
The proof of \Cref{prop:justif-model-xi-ineq} is still valid when considering \(B_{k,j} + \Theta_{k,j}\) instead of \(B_{k,j} + X_{k,j}^{-1}Z_{k,j}\), so that~\eqref{eq:cauchy-decrease} also holds.
As a consequence, \Cref{asm:cauchy-decrease} still holds with \(\varphi\) defined in~\eqref{eq:pd-model-theta}.

Now, our goal is to find a sequence \(\{w_k\} \rightarrow 0\) such that \(w_k \in \nabla f(\bar x_k) - \bar z_k + \partial \psi(\bar s_{k,1}; \bar x_k)\).
Under some additional assumptions on \(\psi\), this will allow us to establish that \Cref{alg:bar-outer} generates iterates that satisfy asymptotically~\eqref{eq:nlp-stationarity-z}.
We begin with preliminary lemmas.

\begin{shadylemma}
  \label{lem:approx-stat-kj}
  Assume that \(s_{k,j,1}\) is not on the boundary of \(\Delta_{k,j} \B \cap \R_{\delta_k}^n(x_{k,j})\) and that \Cref{asm:model} holds.
  Then,
  \begin{equation}
    \label{eq:approx-stat-kj}
    -\nu_{k,j}^{-1} s_{k,j,1} + \mu_k X_{k,j}^{-1}e \in \nabla f(x_{k,j}) + \partial \psi(s_{k,j,1}; x_{k,j}).
  \end{equation} 
\end{shadylemma}

\begin{proof}
  The first step of the proximal gradient method \(s_{k,j,1}\) satisfies~\eqref{eq:skj-1step-pg}.
  According to~\eqref{eq:prox-first-stat}, its first-order optimality conditions are
  \begin{equation}
    \label{eq:skj1-pg-barrier-model}
    0 \in s_{k,j,1} + \nu_{k,j}\nabla \varphi_{\textup{cp}}(0; x_{k,j}) + \partial (\nu_{k,j} \psi(\cdot; x_{k,j}) + \chi(\cdot \mid \Delta_{k,j} \B \cap \R_{\delta_k}^n(x_{k,j})))(s_{k,j,1}).
  \end{equation} 
  As \(s_{k,j,1}\) is not on the boundary of \(\Delta_{k,j} \B \cap \R_{\delta_k}^n(x_{k,j})\), there exists \(r > 0\) such that for all \(s\) in an open ball of center \(s_{k,j,1}\) and radius \(r\), \(\chi(s \mid \Delta_{k,j} \B \cap \R_{\delta_k}^n(x_{k,j})) = 0\).
  Thus, the definition of the subdifferential guarantees that
  \begin{multline*}
    \partial (\nu_{k,j} \psi(\cdot; x_{k,j}) + \chi(\cdot \mid \Delta_{k,j} \B \cap \R_{\delta_k}^n(x_{k,j})))(s_{k,j,1}) = \\
      \nu_{k,j} \partial \psi(s_{k,j,1}; x_{k,j}) + \partial \chi(s_{k,j,1} \mid \Delta_{k,j} \B \cap \R_{\delta_k}^n(x_{k,j})).
  \end{multline*}
  We know that \(\partial \chi(s_{k,j,1} \mid \Delta_{k,j} \B \cap \R_{\delta_k}^n(x_{k,j})) = N_{\Delta_{k,j} \B \cap \R_{\delta_k}^n(x_{k,j})}(s_{k,j,1})\) using \Cref{lem:rn-delta-x-convex}, and \(N_{\Delta_{k,j} \B \cap \R_{\delta_k}^n(x_{k,j})}(s_{k,j,1}) = \{0\}\) because \(s_{k,j,1}\) is not on the boundary of \(\Delta_{k,j} \B \cap \R_{\delta_k}^n(x_{k,j})\).
  Therefore,~\eqref{eq:skj1-pg-barrier-model} simplifies to
  \begin{equation*}
    0 \in \nu_{k,j}^{-1} s_{k,j,1} + \nabla f(x_{k,j}) - \mu_k X_{k,j}^{-1}e + \partial \psi(s_{k,j,1}; x_{k,j})
  \end{equation*}
  using \Cref{asm:model} for \(\nabla \varphi_{\textup{cp}}(0; x_{k,j})\).
\end{proof}

\begin{shadylemma}
  \label{lem:skjk1-epsilondk-ineq}
  Let \Cref{asm:liminf-delta-nu} be satisfied.
  Then, for all \(k \in \N\),
  \begin{equation*}
    \|\bar s_{k,1}\| \le \sqrt{2} \bar \nu_{k} \epsilon_{d,k}.
  \end{equation*}
\end{shadylemma}
\begin{proof}
  Since~\eqref{eq:stat-inner-subproblem} holds,~\eqref{eq:sqrt-nu-xicp-ineq-skj1} leads to
  \begin{equation*}
    \epsilon_{d,k} \ge \bar \nu_k^{-1/2}\xi_{\textup{cp}}(\bar \Delta_k; \bar x_k, \bar \nu_k)^{1/2} \ge \tfrac{1}{\sqrt{2}} \bar \nu_k^{-1}\|\bar s_{k,1}\|,
  \end{equation*}
  which completes the proof.
\end{proof}

The following assumption will be useful to establish that \(\bar s_{k,1}\) converges to zero sufficiently fast to guarantee the convergence of the outer iterations.

\begin{parameterassumption}
  \label{asm:kappa-mdb-eps-cv}
  The sequences \(\{\epsilon_{d,k}\}\) used in \Cref{alg:bar-outer} and \(\{\kappa_{\textup{mdb},k}\}\) from \Cref{prop:min-xkji-kappa-mdb} satisfy
  \begin{equation}
    \kappa_{\textup{mdb},k}^{-1} \epsilon_{d,k} \rightarrow 0.
  \end{equation}
\end{parameterassumption}

To justify that \Cref{asm:kappa-mdb-eps-cv} is reasonable, assume for simplicity that \(i\) is an index such that \((\bar x_k)_i = \kappa_{\textup{mdb}, k}\).
We have
\begin{equation*}
  -e_i^T \bar X_k \bar z_k + \mu_k \le \|\bar X_k \bar z_k - \mu_k e\| \le \epsilon_{p,k},
\end{equation*}
so that 
\begin{equation*}
  (\bar x_k)_i (\bar z_k)_i = e_i^T \bar X_k \bar z_k \ge \mu_k - \epsilon_{p,k},
\end{equation*}
and, if \(\epsilon_{p,k} < \mu_k\),
\begin{equation*}
  \kappa_{\textup{mdb},k}^{-1} = \frac{1}{(\bar x_k)_i} \le \frac{(\bar z_k)_i}{\mu_k - \epsilon_{p,k}}. 
\end{equation*}
We multiply the above inequality by \(\epsilon_{p,k}\) to obtain
\begin{equation*}
  \kappa_{\textup{mdb},k}^{-1} \epsilon_{p,k} \le \frac{(\bar z_k)_i}{\mu_k \epsilon_{p,k}^{-1} - 1}. 
\end{equation*}
If \(\frac{\epsilon_{p,k}}{\mu_k} \rightarrow 0\), e.g., \(\epsilon_{p,k} = \mu_k^{1+\gamma_k}\) with \(0 < \gamma_k < 1\), \(\kappa_{\textup{mdb},k}^{-1} \epsilon_{p,k} \rightarrow 0\).
Then, the choice \(\epsilon_{d,k} = O(\epsilon_{p,k})\) guarantees that \Cref{asm:kappa-mdb-eps-cv} is satisfied.

\begin{shadylemma}
  \label{lem:skj1-not-boundary}
  Let \Cref{asm:model}, \Cref{asm:liminf-delta-nu} and \Cref{asm:kappa-mdb-eps-cv} be satisfied, and for all \(j\), \((f + h)(x_{k, j}) \ge (f + h)_{\textup{low}, k}\).
  Then, there exists \(N \in \mathcal{N}_{\infty}\) such that for all \(k \in N\), \(\bar s_{k,1}\) is not on the boundary of \(\bar \Delta_k \B \cap \R_{\delta_k}^n(\bar x_k)\).
\end{shadylemma}
\begin{proof}
  For all \(s \in \R^n\), \(\|s\|_{\infty} \le \|s\|\), thus \Cref{lem:skjk1-epsilondk-ineq} leads to \(\|\bar s_{k,1}\|_{\infty} \le \|\bar s_{k,1}\| \le \sqrt{2} \bar \nu_k \epsilon_{d,k}\).
  As a consequence, if \(\sqrt{2}\epsilon_{d,k} \bar \nu_k < \Delta_{\min,k}\), then \(\bar s_{k,1}\) is not on the boundary of \(\bar \Delta_k \B\).
  This is certainly true for \(k\) sufficiently large, because \(\bar \Delta > 0\) in \Cref{asm:liminf-delta-nu} and \(\epsilon_{d,k} \rightarrow 0\).

  Now, we show that \(\bar s_{k,1}\) is not on the boundary of \(\R_{\delta_k}^n(\bar x_k)\) if \(k\) is large enough.
  First, we point out that \(\bar \nu_k \not \rightarrow \infty\), because
  \begin{equation*}
    \nu_{k,j} = \frac{1}{\|B_{k,j}\| + \|\Theta_{k,j}\| +  \alpha^{-1}\Delta_{k,j}^{-1}} \le \alpha\Delta_{k,j}.
  \end{equation*}
  Then, we have
  \begin{align*}
    0 \le \frac{|\min_i (\bar s_{k,1})_i|}{\min_i (\bar x_k)_i}
    \le \frac{\|\bar s_{k,1}\|_{\infty}}{\min_i (\bar x_k)_i}
    &\le \frac{\|\bar s_{k,1}\|}{\min_i (\bar x_k)_i} \\
    &\le \frac{\sqrt{2} \bar \nu_k \epsilon_{d,k}}{\min_i (\bar x_k)_i} && \text{using \Cref{lem:skjk1-epsilondk-ineq}}\\
    &\le \frac{\sqrt{2} \bar \nu_k \epsilon_{d,k}}{\kappa_{\textup{mdb},k}} && \text{using \Cref{prop:min-xkji-kappa-mdb}}, 
  \end{align*}
  and \Cref{asm:kappa-mdb-eps-cv} indicates that \(\sqrt{2} \kappa_{\textup{mdb},k}^{-1} \bar \nu_k \epsilon_{d,k} \to 0\).
  As \(\delta_k \rightarrow \bar \delta < 1\), the inequality
  \begin{equation}
    \label{eq:ineq-min-xskj}
    \frac{\min_i (\bar s_{k,1})_i}{\min_i (\bar x_k)_i} + 1 > \delta_k
  \end{equation}
  is satisfied if \(k\) is large enough, and
  \begin{equation}
    \begin{aligned}
      \min_i (\bar s_{k,1} + \bar x_k)_i &\ge \min_i (\bar s_{k,1})_i + \min_i (\bar x_k) & &\text{ by properties of the \(\min\)}\\
      &> \delta_k \min_i (\bar x_k) & &\text{ with~\eqref{eq:ineq-min-xskj}}.\\
    \end{aligned}
  \end{equation}
  Therefore, \(\bar s_{k,1}\) is not on the boundary of \(\R_{\delta_k}(\bar x_k)\) if \(k\) is large enough.
  We conclude that there exists \(N \in \mathcal{N}_{\infty}\) such that for all \(k \in N\), \(\bar s_{k,1}\) is not on the boundary of \(\bar \Delta_k \B \cap \R_{\delta_k}(\bar x_k)\).
\end{proof}

\begin{shadytheorem}
  \label{thm:cv-outer-wk}
  Let \Cref{asm:model}, \Cref{asm:liminf-delta-nu} and \Cref{asm:kappa-mdb-eps-cv} be satisfied, and for all \(j\), \((f + h)(x_{k, j}) \ge (f + h)_{\textup{low}, k}\).
  We define
  \begin{equation}
    \label{eq:def-wk}
    w_k := -\bar \nu_k^{-1} \bar s_{k,1} + \mu_k \bar X_k^{-1}e - \bar z_k.
  \end{equation}
  Then, there exists a subsequence \(N \in \mathcal{N}_{\infty}\) such that for all \(k \in N\),
  \begin{equation}
    \label{eq:asympt-stat-wkj}
    w_k \in \nabla f(\bar x_k) - \bar z_k + \partial \psi(\bar s_{k,1}; \bar x_k),
  \end{equation}
  and
  \begin{equation}
    \label{eq:nkjk-skjk1-lim}
    \|w_k\| \le \sqrt{2}\epsilon_{d,k} + \kappa_{\textup{mdb},k}^{-1}\epsilon_{p,k} \rightarrow 0.
  \end{equation}
\end{shadytheorem}
\begin{proof}
  \Cref{lem:skj1-not-boundary} indicates that there is a subsequence \(N \in \mathcal{N}_{\infty}\) such that, for \(k \in N\), \(\bar s_{k,1}\) is not on the boundary of \(\bar \Delta_k \B \cap \R_{\delta_k}^n(\bar x_k)\).
  Thus, \Cref{lem:approx-stat-kj} holds, and
  \begin{equation}
    \label{eq:partia-psi-approx-sol}
    -\bar \nu_k^{-1} \bar s_{k,1} + \mu_k \bar X_k^{-1}e - \bar z_k \in \nabla f(\bar x_k) - \bar z_k + \partial \psi(\bar s_{k,1}; \bar x_k).
  \end{equation}
  With \(w_k\) defined in~\eqref{eq:def-wk}, we have~\eqref{eq:asympt-stat-wkj}.
  Now, for all \(i \in \{1, \ldots, n\}\), we use \Cref{prop:min-xkji-kappa-mdb} to establish that
  \begin{equation*}
    (\mu_k \bar X_k^{-1}e - \bar z_k)_i = (\mu_k  - (\bar x_k)_i (\bar z_k)_i) / (\bar x_k)_i \le (\mu_k  - (\bar x_k)_i (\bar z_k)_i) / \kappa_{\textup{mdb},k},
  \end{equation*}
  and, by summing the square of the above inequality for all \(i \in \{1, \ldots, n\}\),
  \begin{equation}
    \label{eq:compl-square-ineq-eps}
    \begin{aligned}
      \|\mu_k \bar X_k^{-1}e - \bar z_k\|^2 &= \sum_{i=1}^n (\mu_k \bar X_k^{-1}e - \bar z_k)_i^2 \\
      &\le \sum_{i=1}^n (\mu_k  - (\bar x_k)_i (\bar z_k)_i)^2 / \kappa_{\textup{mdb},k}^2 \\
      &= \|\mu_k e - \bar X_k \bar z_k\|^2 / \kappa_{\textup{mdb},k}^2 \\
      &\le \kappa_{\textup{mdb},k}^{-2} \epsilon_{p,k}^2 \text{ because~\eqref{eq:central-path-inner} holds},
    \end{aligned}
  \end{equation}
  so that \(\|\mu_k \bar X_k^{-1}e - \bar z_k\| \le \kappa_{\textup{mdb},k}^{-1} \epsilon_{p,k}\).
  As \(\epsilon_{p,k} \rightarrow 0\) and \Cref{asm:kappa-mdb-eps-cv} holds, we deduce that \(\|\mu_k \bar X_k^{-1}e - \bar z_k\| \rightarrow 0\).

  Finally,
  \begin{align*}
    \|w_k\| &= \|- \bar \nu_k^{-1} \bar s_{k,1} + \mu_k \bar X_k^{-1}e - \bar z_k\| \\
    &\le \bar \nu_k^{-1} \|\bar s_{k,1}\| + \|\mu_k \bar X_k^{-1}e - \bar z_k\| \\
    &\le \sqrt{2} \epsilon_{d,k} + \kappa_{\textup{mdb},k}^{-1} \epsilon_{p,k} \underset{k \rightarrow +\infty}{\longrightarrow} 0,
  \end{align*}
  where we used \Cref{lem:skjk1-epsilondk-ineq} and~\eqref{eq:compl-square-ineq-eps} in the last inequality.
\end{proof}

Now, we present two assumptions on \(\psi\).
The first will not be necessary for the remaining results of this subsection, except as one of the justifications for the second assumption.
However, it will be used in \Cref{sec:new-crit-meas}.

\begin{modelassumption}
  \label{asm:psi-epi-continuity}
  \(x \mapsto \psi(\cdot; x)\) is epi-continuous on \(\R_+^n\).
\end{modelassumption}

\Cref{asm:psi-epi-continuity} holds if \(\psi\) is continuous on \(\R_+^n \times \R_+^n\), but this condition is only sufficient, not necessary \citep[Exercise~\(7.40\)]{rtrw}.
Let us consider the case where \(\psi(\cdot; x) = s \mapsto h(x + s)\).
Because \(h\) is lsc, its epigraph is closed, thus the sequence of functions \(\{h, h, \ldots\}\) satisfies \(\{h, h, \ldots\} \eto h\).
Let \(\bar x \in \R_+^n\) and \(\{x_k\} \rightarrow \bar x\).
\citep[Exercise~\(7.8d\)]{rtrw} indicates that for \(h_k = s \mapsto h(x_k + s)\) and \(\bar h = s \mapsto h(\bar x + s)\), \(h_k \eto \bar h\).
Since the latter is true for all \(\bar x \in \R_+^n\), we conclude that \Cref{asm:psi-epi-continuity} is satisfied.

\begin{modelassumption}
  \label{asm:subgrad-cv}
  For any sequences \(\{s_k\} \rightarrow 0\) and \(\{x_k\} \rightarrow \bar x \ge 0\) such that \(x_k + s_k > 0\) for all \(k \in \N\),
  \begin{equation}
    \limsup_{k \rightarrow \infty} \partial \psi(s_k; x_k) \subset \partial \psi(0; \bar x) = \partial h(\bar x).
  \end{equation}
\end{modelassumption}
We present some cases for which \Cref{asm:subgrad-cv} holds.
\begin{itemize}
  \item When \(\glim_{k \rightarrow \infty} \partial \psi(\cdot; x_k) = \partial \psi(\cdot; \bar x)\).
  Attouch's theorem \citep{attouch-1977} (also written in \citep[Theorem~\(12.35\)]{rtrw}) indicates that this condition is satisfied when \(\psi(\cdot; x_k)\) and \(\psi(\cdot; \bar x)\) are proper, lsc, convex functions with \(\psi(\cdot; x_k) \eto \psi(\cdot; \bar x)\) (i.e., \Cref{asm:psi-epi-continuity} holds).
  An extension to non-convex functions under some more sophisticated assumptions is established by \citet{poliquin-1992}.
  \item When \(\psi(s; x) = h(x + s)\) and \(h(x_k + s_k) \rightarrow h(\bar x)\) (e.g., \(h = \|\cdot\|_1\)), using \Cref{prop:subgradient-osc} applied to \(\{x_k + s_k\} \to \bar x\).
  \item When \(\psi(s; x) = h(x + s)\) but \(h\) is not continuous, we may still be able to show that \Cref{asm:subgrad-cv} holds.
  For example, with \(h = \|\cdot\|_0\), \citep[Theorem~\(1\)]{le-2013} shows that \(\partial h(x) = \{v \mid v_i = 0 \text{ if } x_i \neq 0\}\).
  Thus, \(\partial \psi(s_k; x_k) = \partial (\|\cdot\|_0)(x_k + s_k) = \{0\} \subset \partial h(\bar x)\).
\end{itemize}

\begin{shadycorollary}
  \label{cor:stationarity}
  Under the assumptions of \Cref{thm:cv-outer-wk} and \Cref{asm:subgrad-cv}, let \(\bar x\) and \(\bar z\) be limit points of \(\{\bar x_k\}\) and \(\{\bar z_k\}\), respectively.
  Then,
  \begin{equation}
    0 \in \nabla f(\bar x) - \bar z + \partial h(\bar x) \quad \text{and} \quad \bar X \bar Z e = 0.
  \end{equation}
  In this case, when the CQ is satisfied at \(\bar x\), \(\bar x\) is first-order stationary for~\eqref{eq:nlp}.
\end{shadycorollary}

\begin{proof}
  In \Cref{thm:cv-outer-wk}, \(N\) can be chosen such that \(\bar x_k \underset{k \in N}{\longrightarrow} \bar x\) and \(\bar z_k \underset{k \in N}{\longrightarrow} \bar z\).
  We apply \Cref{asm:subgrad-cv} to~\eqref{eq:asympt-stat-wkj} and~\eqref{eq:nkjk-skjk1-lim} to deduce
  \begin{equation}
    \label{eq:zero-in-P-bar}
    -\nabla f(\bar x) + \bar z \in \partial \psi(0; \bar x) = \partial h(\bar x),
  \end{equation}
  which indicates that \(0 \in P^{\mathcal{L}}(\bar \Delta; \bar x, \bar z, \bar \nu)\).
  The condition~\eqref{eq:central-path-inner} implies that \(\bar X \bar Z e = 0\).
  When the CQ is satisfied, we can use \Cref{lem:xi-stationarity} to conlude that \(\bar x\) is first-order stationary for~\eqref{eq:nlp}.
\end{proof}

In \Cref{thm:cv-outer-wk}, \(w_k\) evokes of the concept of \((\epsilon_p, \epsilon_d)\)-KKT optimality for interior-point methods introduced in \citep[Definition~\(2.1\)]{demarchi-themelis-2022}.
We slightly modify this concept in the following definition.

\begin{definition}
  \label{def:eps-kkt-stat}
  Let \(\epsilon_p, \epsilon_d \ge 0\), \(x \ge 0\), and \(\Delta, \nu > 0\).
  \(x\) is said to be \((\epsilon_p, \epsilon_d)\)-KKT optimal if there exist \(z \ge 0\) and \(v \in \partial h(x)\) such that
  \begin{equation}
    \label{eq:eps-d-stat}
    \|\nabla f(x) - z + v\| \le \epsilon_d,
  \end{equation}
  and, for all \(i \in \{1, \ldots, n\}\),
  \begin{equation}
    x_i z_i \le \epsilon_p.
  \end{equation}
\end{definition}

The main modification to the original formulation in \citep[Definition~\(2.1\)]{demarchi-themelis-2022} is that we require \(x_i z_i \le \epsilon_p\) instead of \(\min (x_i, z_i) \le \epsilon_p\) for all \(i\), but this is linked to our different choices of stopping condition for the complementary slackness.
The first part of the definition~\eqref{eq:eps-d-stat} is similar to the \(\epsilon\)-stationarity \citep[Definition~\(4.5\)]{demarchi-themelis-2022b} for more general problems.

\Cref{asm:model} does not necessarily guarantee that \(\partial \psi(s_{k,j,1}; x_{k,j}) = \partial h(x_{k,j} + s_{k,j,1})\).
Thus, for \(k \in N\) where \(N\) is a subsequence introduced in \Cref{thm:cv-outer-wk}, we cannot use \Cref{thm:cv-outer-wk} to measure the \((\epsilon_p, \epsilon_d)\)-KKT optimality of \(\bar x_k\).
Let
\begin{equation}
  \label{eq:def-eps-hk}
  \epsilon_{h,k} = \dist(w_k - \nabla f(\bar x_k) + \bar z_k, \, \partial h(\bar x_k + \bar s_{k,1})).
\end{equation}
As \Cref{thm:cv-outer-wk} indicates that \(w_k - f(\bar x_k) + \bar z_k \in \partial \psi(\bar s_{k,1}; \bar x_k)\), we can obtain a measure of \((\epsilon_p, \epsilon_d)\)-KKT optimality which depends on \(\epsilon_{h,k}\).
When all the elements of \(\partial \psi(\bar s_{k,1}; \bar x_k)\) are close to an element of \(\partial h(\bar x_k + \bar s_{k,1})\), we expect \(\epsilon_{h,k}\) to be small.
In particular, if \(\partial \psi(\bar s_{k,1}; \bar x_k) \subseteq \partial h(\bar x_k + \bar s_{k,1})\), \(\epsilon_{h,k} = 0\).
\begin{shadytheorem}
  \label{thm:eps-kkt-stat}
  Let the assumptions of \Cref{thm:cv-outer-wk} be satisfied, and \(\epsilon_{h,k}\) be defined in~\eqref{eq:def-eps-hk}.
  Then, there exists \(N \in \mathcal{N}_{\infty}\) such that for all \(k \in N\), \(\bar x_k + \bar s_{k,1}\) is \((\bar \epsilon_{p,k}, \bar \epsilon_{d,k})\)-KKT optimal with constants
  \begin{align*}
    \bar \epsilon_{p,k} &= \epsilon_{p,k} + \sqrt{n}\mu_k + \sqrt{2}\bar \nu_k\epsilon_{d,k}\|\bar z_k\| \\
    \bar \epsilon_{d,k} &= \epsilon_{h,k} + \sqrt{2}\epsilon_{d,k}(1 + \bar \nu_k L_f) + \kappa_{\textup{mdb},k}^{-1}\epsilon_{p,k}.
  \end{align*}
\end{shadytheorem}

\begin{proof}
  \Cref{thm:cv-outer-wk} guarantees that, for all \(k\) in a subsequence \(N \in \mathcal{N}_{\infty}\),~\eqref{eq:asympt-stat-wkj} holds, i.e.
  \begin{equation*}
    w_k - \nabla f(\bar x_k) + \bar z_k \in \partial \psi(\bar s_{k,1}; \bar x_k),
  \end{equation*}
  where \(w_k\) is defined in~\eqref{eq:def-wk}.
  As \(\partial h(\bar x_k + \bar s_{k,1})\) is closed, we can choose \(v_k \in \partial h(\bar x_k + \bar s_{k,1})\) such that, for \(y_k := v_k - w_k + \nabla f(\bar x_k) - \bar z_k\), we have \(\|y_k\| = \epsilon_{h,k}\).
  Then,
  \begin{equation*}
    (v_k - w_k + \nabla f(\bar x_k) - \bar z_k) + w_k  \in \nabla f(\bar x_k) - \bar z_k + \partial h(\bar x_k + \bar s_{k,1}),
  \end{equation*}
  which we may rewrite as
  \begin{equation*}
    y_k + w_k + \nabla f(\bar x_k + \bar s_{k,1}) - \nabla f(\bar x_k) \in \nabla f(\bar x_k + \bar s_{k,1}) - \bar z_k + \partial h(\bar x_k + \bar s_{k,1}),
  \end{equation*}
  and
  \begin{equation*}
      \|y_k + w_k + \nabla f(\bar x_k + \bar s_{k,1}) - \nabla f(\bar x_k)\| \le \epsilon_{h,k} + \|w_k\| + L_f \|\bar s_{k,1}\|.
  \end{equation*}
  \Cref{lem:skjk1-epsilondk-ineq} implies that \(L_f \|\bar s_{k,1}\| \le \sqrt{2} L_f \bar \nu_k \epsilon_{d,k}\).
  We combine the latter inequality with~\eqref{eq:nkjk-skjk1-lim} in \Cref{thm:cv-outer-wk} to obtain
  \begin{equation}
    \label{eq:eps-d-ineq-kkt-stat}
    \|y_k + w_k + \nabla f(\bar x_k + \bar s_{k,1}) - \nabla f(\bar x_k)\| \le \epsilon_{h,k} + \sqrt{2}\epsilon_{d,k} (1 + \bar \nu_k L_f) + \kappa_{\textup{mdb},k}^{-1}\epsilon_{p,k}.
  \end{equation}
  Now, for all \(i \in \{1, \ldots, n\}\),
  \begin{equation}
    \label{eq:eps-p-ineq-kkt-stat}
    \begin{aligned}
      (\bar x_k + \bar s_{k,1})_i (\bar z_k)_i &\le \|(\bar X_k + \bar S_{k,1}) \bar z_k\| \\
      &\le \|\bar X_k \bar z_k - \mu_k e\| + \|\mu_k e\| + \|\bar S_{k,1}\bar z_k\|\\
      &\le \epsilon_{p,k} + \sqrt{n} \mu_k + \|\bar s_{k,1}\|\|\bar z_k\| \\
      &\le \epsilon_{p,k} + \sqrt{n} \mu_k + \sqrt{2}\bar \nu_k\epsilon_{d,k}\|\bar z_k\|.
    \end{aligned}
  \end{equation}
  We use~\eqref{eq:eps-d-ineq-kkt-stat},~\eqref{eq:eps-p-ineq-kkt-stat} and \Cref{def:eps-kkt-stat} to conclude.
\end{proof}

If \(\bar z_k\) is bounded, \(\bar \epsilon_{p,k} \underset{N}{\longrightarrow} 0\) in \Cref{thm:eps-kkt-stat}.
If \(\epsilon_{h,k} \rightarrow 0\), we also have \(\bar \epsilon_{d,k} \underset{N}{\longrightarrow} 0\).

\subsection{Convergence with a new criticality measure}
\label{sec:new-crit-meas}

Now, instead of using \(\nu_{k,j}^{-1/2} \xi_{\textup{cp}}(\Delta_{k,j}; x_{k,j}, \nu_{k,j})^{1/2}\) (involving \(\varphi_{\textup{cp}}(\cdot; x_{k,j})\)) for the criticality measure of \Cref{alg:bar-inner}, we would like to use a measure based upon \(\varphi^{\mathcal{L}} (\cdot; x_{k,j}, z_{k,j})\) defined in~\eqref{eq:bar-varphi}.
The reason behind this choice is inspired from the criticality measure \(\|\nabla f(x_{k,j}) - z_{k,j}\|\) used in primal-dual trust region algorithms in the smooth case, instead of \(\|\nabla f(x_{k,j}) - \mu_k X_{k,j}^{-1}e\|\) used in primal algorithms, see for example \citep[Algorithm~\(13.6.2\)]{conn-gould-toint-2000}. 
We may expect that this choice results in fewer iterations of \Cref{alg:bar-inner} when \(x_{k,j}\) and \(z_{k,j}\) are close to a solution of~\eqref{eq:nlp}, because (we express this idea with smooth notations for now), if \(j_k\) is the index for which the stopping criteria of \Cref{alg:bar-inner} are met, \(\|\nabla f(\bar x_k) - \bar z_k\| = \|\nabla f(x_{k+1, 0}) - z_{k+1,0}\|\), whereas \(\|\nabla f(\bar x_k) - \mu_k \bar X_k^{-1}e\| \neq \|\nabla f(x_{k+1,0}) - \mu_{k+1} X_{k+1,0}^{-1}e\|\).
However, to change the stopping criterion, we need the following convexity assumption.

\begin{modelassumption}
  \label{asm:psi-convex}
  For a sequence \(\{x_{k,j}\}_{j}\) generated by \Cref{alg:bar-inner} at iteration \(k\), \(\psi(\cdot; x_{k,j})\) is convex for all \(j\).
\end{modelassumption}

If \(\psi(s; x) = h(x + s)\) and \(h\) is convex, \Cref{asm:psi-convex} holds.

In this section, we define
\begin{equation}
  \label{eq:bar-skj-L}
  s_{k,j}^{\mathcal{L}} \in \argmin{s} m^{\mathcal{L}}(s; x_{k,j}, z_{k,j}, \nu_{k,j}) + \chi(s \mid \Delta_{k,j} \B \cap \R_{\delta_k}^n(x_{k,j})),
\end{equation}
where \(m^{\mathcal{L}}(s; x_{k,j}, z_{k,j} \nu_{k,j})\) is defined in~\eqref{eq:bar-model}, and
\begin{equation}
  \label{eq:def-xi-deltak}
  \xi_{\delta_k}^{\mathcal{L}}(\Delta_{k,j}; x_{k,j}, z_{k,j}, \nu_{k,j}) = (f + \phi_k + h)(x_{k,j}) - (\varphi^{\mathcal{L}}(s_{k,j}^{\mathcal{L}}; x_{k,j}, z_{k,j}) - \psi(s_{k,j}^{\mathcal{L}}; x_{k,j})),
\end{equation}
where \(\varphi^{\mathcal{L}}(\cdot; x_{k,j}, z_{k,j})\) is defined in~\eqref{eq:bar-varphi}.
We point out that \(\xi_{\delta_k}^{\mathcal{L}}\) and \(\xi^{\mathcal{L}}\) defined in~\eqref{eq:def-xibar} are almost identical, the latter being computed by replacing \(\chi(s \mid \Delta_{k,j} \B \cap \R_{\delta_k}^n(x_{k,j}))\) by \(\chi(s \mid \Delta_{k,j} \B) + \chi(x_{k,j} + s \mid \R_+^n)\) in~\eqref{eq:bar-skj-L}.
\Cref{asm:psi-convex} will be useful in \Cref{thm:inner-cv-xi-deltak}, which is crucial for our analysis of \Cref{alg:bar-outer-z} because it establishes the convergence of the inner iterations with \(\xi_{\delta_k}^{\mathcal{L}}\).

\Cref{alg:bar-outer-z} resembles \Cref{alg:bar-outer}, except for the stopping criterion~\eqref{eq:stat-inner-subproblem}.
Our ultimate goal in this subsection is to show that replacing~\eqref{eq:stat-inner-subproblem} in \Cref{alg:bar-outer} by~\eqref{eq:stat-inner-subproblem-z} in \Cref{alg:bar-outer-z} maintains similar convergence properties to those of \Cref{sec:convergence-outer}, but, as illustrated in \Cref{sec:numerical}, performs better in practice.

\begin{algorithm}[h]
  \caption[caption]{%
    Nonsmooth interior-point method (outer iteration) with stopping criteria based upon \(\xi_{\delta_k}^{\mathcal{L}}\) in~\eqref{eq:def-xi-deltak}.%
    \label{alg:bar-outer-z}
  }
  \begin{algorithmic}[1]
    \State Choose \(\epsilon > 0\), sequences \(\{\mu_k\} \searrow 0\), \(\{\epsilon_{d,k}\} \searrow 0\), \(\{\epsilon_{p,k}\} \searrow 0\), and \(\{\delta_k\} \rightarrow \bar \delta \in [0,1)\) with \(\delta_k \in (0,1)\) for all \(k\).
    \State Choose \(x_{0, 0} \in \R^n_{++}\) where \(h\) is finite.
    \For{\(k = 0, 1, \ldots\)}
      \State Compute an approximate solution \(x_k := x_{k,j}\) to~\eqref{eq:bar-subproblem} and \(z_k := z_{k,j}\) in the sense that
      \begin{equation}
        \label{eq:stat-inner-subproblem-z}
        \nu_{k,j}^{-1/2} \xi_{\delta_k}^{\mathcal{L}}(\Delta_{k,j}; x_{k,j}, \nu_{k,j})^{1/2} \leq \epsilon_{d,k},
      \end{equation}
      and~\eqref{eq:central-path-inner} holds, that we recall in the following inequality for conveniency
      \begin{equation*}
        \|\mu_k e - X_{k,j} z_{k,j}\| \le \epsilon_{p,k}.
      \end{equation*}
      \State Set \(x_{k+1, 0} := x_k\).
    \EndFor
  \end{algorithmic}
\end{algorithm}

The following result shows that the inner iteration terminates finitely.

\begin{shadytheorem}
  \label{thm:inner-cv-xi-deltak}
  Under the assumptions of \Cref{prop:complementarity}, \Cref{asm:psi-epi-continuity} and \Cref{asm:psi-convex}, if \(\{x_{k,j}\}_j\) possesses a limit point \(x_k^*\), then
  \begin{equation}
    \label{eq:xibar-cv-0-j}
    \liminf_{j \rightarrow +\infty} \xi_{\delta_k}^{\mathcal{L}}(\Delta_{k,j}; x_{k,j}, z_{k,j}, \nu_{k,j}) = 0,
  \end{equation}
  and
  \begin{equation}
    \label{eq:skjL-cv-0-j}
    \liminf_{j \rightarrow +\infty} \|s_{k,j}^{\mathcal{L}}\| = 0,
  \end{equation}
  where \(s_{k,j}^{\mathcal{L}}\) is defined in~\eqref{eq:bar-skj-L}.
\end{shadytheorem}

\begin{proof}
  As \(\Delta_{k,j} \in [\Delta_{\min, k}, \Delta_{\max}]\) and \(\nu_{k,j} \in [\nu_{\min,k}, 1]\), there exists an infinite subsequence \(N\) such that \(\Delta_{k,j} \underset{j \in N}{\longrightarrow} \Delta_k^*\), \(\nu_{k,j} \underset{j \in N}{\longrightarrow} \nu_k^*\), \(x_{k,j} \underset{j \in N}{\longrightarrow} x_k^*\), and with \Cref{prop:complementarity} \(z_{k,j} \underset{j \in N}{\longrightarrow} z_k^*\) with \(X_k^* Z_k^* e = \mu_k e\).
  By continuity of the \(\min\), \(\R_{\delta_k}^n(x_{k,j}) \underset{j \in N}{\longrightarrow} \R_{\delta_k}^n(x_k^*)\).
  The sets \(\Delta_{k,j} \B\) and \(\R_{\delta_k}^n(x_{k,j})\) are convex (using \Cref{lem:rn-delta-x-convex} for the latter).
  Since \(\Delta_k^* \B\) and \(\R_{\delta_k}^n(x_k^*)\) are convex and cannot be separated, we use \citep[Theorem~\(4.33\)]{rtrw} to conclude that
  \begin{equation*}
    \Delta_{k,j} \B \cap \R_{\delta_k}^n(x_{k,j}) \underset{j \in N}{\longrightarrow} \Delta_k^* \B \cap \R_{\delta_k}^n(x_k^*).
  \end{equation*}
  With \citep[Theorem~\(7.4f\)]{rtrw}, we deduce
  \begin{equation*}
    \elim_{j \in N} \chi(\cdot \mid \Delta_{k,j} \B \cap \R_{\delta_k}^n(x_{k,j})) = \chi(\cdot \mid \Delta_k^* \B \cap \R_{\delta_k}^n(x_k^*)).
  \end{equation*}
  Thanks to \Cref{prop:complementarity} and the smoothness of \(\varphi_{\textup{cp}}(\cdot; x)\) and \(\varphi^{\mathcal{L}}(\cdot; x, z)\), we also have
  \begin{equation}
    \label{eq:elim-varphi-j}
    \elim_{j \in N} \varphi_{\textup{cp}}(\cdot; x_{k,j}) = \elim_{j \in N} \varphi^{\mathcal{L}}(\cdot; x_{k,j}, z_{k,j}) = \varphi^{\mathcal{L}}(\cdot; x_k^*, z_k^*).
  \end{equation}
  The functions \(\varphi_{\textup{cp}}(\cdot; x_{k,j})\), \(\varphi^{\mathcal{L}} (\cdot; x_{k,j}, z_{k,j})\) and \(\varphi^{\mathcal{L}} (\cdot; x_k^*, z_k^*)\) are all convex because they are linear.
  \Cref{asm:psi-epi-continuity} implies that \(\elim_{j\in N} \psi(\cdot; x_{k,j}) = \psi(\cdot; x_k^*)\).
  \Cref{asm:psi-convex} and \citep[Theorem~\(7.46\)]{rtrw} lead to
  \begin{equation}
    \label{eq:elim-psi-chi}
    \elim_{j \in N} \psi(\cdot; x_{k,j}) + \chi(\cdot \mid \Delta_{k,j} \B \cap \R_{\delta_k}^n(x_{k,j})) = \psi(\cdot; x_k^*) + \chi(\cdot \mid \Delta_k^* \B \cap \R_{\delta_k}^n(x_k^*)),
  \end{equation}
  and the above functions are all convex.
  We deduce from~\eqref{eq:elim-varphi-j},~\eqref{eq:elim-psi-chi} and \citep[Theorem~\(7.46\)]{rtrw} that
  \begin{equation*}
    \elim_{j \in N} m_{\textup{cp}}(\cdot; x_{k,j}, \nu_{k,j}) = \elim_{j \in N} m^{\mathcal{L}}(\cdot; x_{k,j}, z_{k,j}, \nu_{k,j}) = m^{\mathcal{L}}(\cdot; x_k^*, z_k^*, \nu_k^*),
  \end{equation*}
  where \(m^{\mathcal{L}}\) is defined in~\eqref{eq:bar-model}.
  The sequences \(m_{\textup{cp}}(\cdot; x_{k,j}, \nu_{k,j}) + \chi(\cdot \mid \Delta_{k,j} \B \cap \R_{\delta_k}^n(x_{k,j}))\) and \(m^{\mathcal{L}}(\cdot; x_k^*, z_k^*, \nu_k^*) + \chi(\cdot \mid \Delta_{k,j} \B \cap \R_{\delta_k}^n(x_{k,j}))\) are level-bounded because of the indicators.
  As in~\eqref{eq:skj-1step-pg}, we have
  \begin{multline*}
    \prox{\nu_k^* \psi(\cdot; x_k^*) + \chi(\cdot \mid \Delta_k^* \B \cap \R_{\delta_k}^n(x_k^*))}(-\nu_k^*\nabla \varphi_{\textup{cp}}(0; x_k^*, \nu_k^*)) = \\
    \arg \min_s m^{\mathcal{L}}(s; x_k^*, z_k^*, \nu_k^*) + \chi(s \mid \Delta_k^* \B \cap \R_{\delta_k}^n(x_k^*)),
  \end{multline*}
  and the above problem is single valued because of \citep[Theorem~\(2.26a\)]{rtrw}.
  Let \(s^*\) denote its only solution.
  \Cref{thm:epi-cv-minimization}, and specifically~\eqref{eq:limsup-argmin-phik}, implies that the sequences
  \begin{equation*}
    s_{k,j,1} \in \arg \min_s m_{\textup{cp}}(s; x_{k,j}, \nu_{k,j}) + \chi(s \mid \Delta_{k,j} \B \cap \R_{\delta_k}^n(x_{k,j})),
  \end{equation*}
  and
  \begin{equation*}
    s_{k,j}^{\mathcal{L}} \in \arg \min_s m^{\mathcal{L}}(s; x_{k,j}, z_{k,j}, \nu_{k,j}) + \chi(s \mid \Delta_{k,j} \B \cap \R_{\delta_k}^n(x_{k,j}))
  \end{equation*}
  have the same limit \(s^*\).
  We have shown in \Cref{lem:skj-0} that \(s_{k,j,1} \underset{j \rightarrow +\infty}{\longrightarrow} 0\).
  Thus, \(s^* = 0\).
  Finally, we have
  \begin{equation*}
    \xi_{\textup{cp}}(\Delta_{k,j}; x_{k,j}, \nu_{k,j}) = (f + h + \phi_k)(x_{k,j}) - m_{\textup{cp}}(s_{k,j,1}; x_{k,j}, \nu_{k,j}) + \tfrac{1}{2}\nu_{k,j}^{-1}\|s_{k,j,1}\|^2.
  \end{equation*}
  As \(\xi_{\textup{cp}}(\Delta_{k,j}; x_{k,j}, \nu_{k,j}) \underset{j \rightarrow +\infty}{\longrightarrow} 0\) by \Cref{prop:cv-inner}, and \(\|s_{k,j,1}\| \underset{j \rightarrow +\infty}{\longrightarrow} 0\), we deduce that
  \begin{equation*}
    m_{\textup{cp}}(s_{k,j,1}; x_{k,j}, \nu_{k,j}) \underset{j \in N}{\longrightarrow} (f + h + \phi_k)(x_k^*).
  \end{equation*}
  Using~\eqref{eq:cv-inf-phik} in \Cref{thm:epi-cv-minimization}, we have
  \begin{equation}
    \label{eq:lim-mcp-mL}
    \lim_{j \in N} m_{\textup{cp}}(s_{k,j,1}; x_{k,j}, \nu_{k,j}) = \lim_{j \in N} m^{\mathcal{L}}(s_{k,j}^\mathcal{L}; x_{k,j}, z_{k,j}, \nu_{k,j}) = (f + h + \phi_k)(x_k^*).
  \end{equation}
  The expression of~\(\xi_{\delta_k}^{\mathcal{L}}\) in~\eqref{eq:def-xi-deltak} can also be written as
  \begin{equation*}
    \xi_{\delta_k}^{\mathcal{L}}(\Delta_{k,j}; x_{k,j}, z_{k,j}, \nu_{k,j}) = (f + h + \phi_k)(x_{k,j}) - m^{\mathcal{L}}(s_{k,j}^{\mathcal{L}}; x_{k,j}, z_{k,j}, \nu_{k,j}) + \tfrac{1}{2}\nu_{k,j}^{-1}\|s_{k,j}^{\mathcal{L}}\|^2.
  \end{equation*}
  By injecting the limit of \(s_{k,j}^{\mathcal{L}}\) and~\eqref{eq:lim-mcp-mL} in the above equation, we obtain~\eqref{eq:xibar-cv-0-j}.
\end{proof}

From this point on, \(j_k\) denotes the number of iterations performed by \Cref{alg:bar-inner} at iteration \(k\) with the inner stopping criteria from \Cref{alg:bar-outer-z}, and we use again the notation \(\bar x_k = x_{k,j_k}\), \(\bar z_k = z_{k,j_k}\), \(\bar s_{k,1} = s_{k,j_k,1}\), \(\bar \Delta_k = \Delta_{k,j_k}\), \(\bar \nu_{k} = \nu_{k,j_k}\), with the addition of \(\bar s_k^{\mathcal{L}} := s_{k,j_k}^{\mathcal{L}}\).
The following three lemmas are analogous to \Cref{lem:approx-stat-kj}, \Cref{lem:skjk1-epsilondk-ineq} and \Cref{lem:skj1-not-boundary}.

\begin{shadylemma}
  \label{lem:approx-stat-kj-z}
  Assume that \(s_{k,j}^{\mathcal{L}}\) is not on the boundary of \(\Delta_{k,j} \B \cap \R_{\delta_k}^n(x_{k,j})\) and that \Cref{asm:model} holds.
  Then,
  \begin{equation}
    \label{eq:approx-stat-kj-z}
    -\nu_{k,j}^{-1} s_{k,j}^{\mathcal{L}} \in \nabla f(x_{k,j}) - z_{k,j} + \partial \psi(s_{k,j}^{\mathcal{L}}; x_{k,j}).
  \end{equation} 
\end{shadylemma}

\begin{proof}
  The first-order stationarity condition of~\eqref{eq:bar-skj-L} is
  \begin{equation*}
    -\nu_{k,j}^{-1} s_{k,j}^{\mathcal{L}} \in \nabla f(x_{k,j}) - z_{k,j} + \partial \psi(s_{k,j}^{\mathcal{L}}; \bar x_k) + \partial \chi(s_{k,j}^{\mathcal{L}} \mid \Delta_{k,j} \B \cap \R_{\delta_k}^n(x_{k,j})).
  \end{equation*}
  The same analysis as in the proof of~\Cref{lem:approx-stat-kj} establishes that
  \begin{equation*}
    \partial \chi(s_{k,j}^{\mathcal{L}} \mid \Delta_{k,j} \B \cap \R_{\delta_k}^n(x_{k,j})) = \{0\},
  \end{equation*}
  so that~\eqref{eq:approx-stat-kj-z} holds.
\end{proof}

\begin{shadylemma}
  \label{lem:bar-skjk1-epsilondk-ineq}
  Let \Cref{asm:liminf-delta-nu} be satisfied.
  Then, for all \(k \in \N\),
  \begin{equation*}
    \|\bar s_k^{\mathcal{L}}\| \le \sqrt{2} \bar \nu_k \epsilon_{d,k}.
  \end{equation*}
\end{shadylemma}
\begin{proof}
  The bound
  \begin{equation}
    \label{eq:xi-deltak-skj-L-ineq}
    \xi_{\delta_k}^{\mathcal{L}}(\Delta_{k,j}; x_{k,j}, z_{k,j}, \nu_{k,j}) \ge \tfrac{1}{2} \nu_{k,j}^{-1}\|s_{k,j}^{\mathcal{L}}\|^2
  \end{equation}
  holds because
  \begin{align*}
    m^{\mathcal{L}}(0; x_{k,j}, z_{k,j}, \nu_{k,j}) &= (f + \phi_k + h)(x_{k,j}) \\
    &\ge m^{\mathcal{L}}(s_{k,j}^{\mathcal{L}}; x_{k,j}, z_{k,j}, \nu_{k,j})\\
    &= (f + \phi_k + h)(x_{k,j}) + (\nabla f(x_{k,j}) + z_{k,j})^T s_{k,j}^{\mathcal{L}} + \tfrac{1}{2}\nu_{k,j}^{-1}\|s_{k,j}^{\mathcal{L}}\|^2.
  \end{align*}
  The stopping criterion~\eqref{eq:stat-inner-subproblem-z} and \eqref{eq:xi-deltak-skj-L-ineq} lead to
  \begin{equation*}
    \tfrac{1}{\sqrt{2}}\bar \nu_k^{-1}\|\bar s_k^{\mathcal{L}}\| \le \bar \nu_k^{-1/2}\xi_{\delta_k}^{\mathcal{L}}(\bar \Delta_k; \bar x_k, \bar z_k, \bar \nu_k)^{1/2} \le \epsilon_{d,k},
  \end{equation*}
  which completes the proof.
\end{proof}

\begin{shadylemma}
  \label{lem:bar-skj1-not-boundary}
  Let \Cref{asm:model}, \Cref{asm:liminf-delta-nu} and \Cref{asm:kappa-mdb-eps-cv} be satisfied, and for all \(j\), \((f + h)(x_{k, j}) \ge (f + h)_{\textup{low}, k}\).
  Then, there exists \(N \in \mathcal{N}_{\infty}\) such that for all \(k \in N\), \(\bar s_k^{\mathcal{L}}\) is not on the boundary of \(\bar \Delta_k \B \cap \R_{\delta_k}^n(\bar x_k)\).
\end{shadylemma}
\begin{proof}
  Since for any \(s \in \R^n\), \(\|s\|_{\infty} \le \|s\|\), \Cref{lem:bar-skjk1-epsilondk-ineq} leads to \(\|\bar s_k^{\mathcal{L}}\|_{\infty} \le \|\bar s_k^{\mathcal{L}}\| \le \sqrt{2} \bar \nu_k \epsilon_{d,k}\).
  Thus, if \(\sqrt{2}\epsilon_{d,k} \bar \nu_k < \Delta_{\min,k}\), \(\bar s_k^{\mathcal{L}}\) is not on the boundary of \(\bar \Delta_k \B\).
  As \(\bar \Delta > 0\) in \Cref{asm:liminf-delta-nu} and \(\epsilon_{d,k} \rightarrow 0\), this is true if \(k\) is sufficiently large.
  The rest of the proof is identical to that of \Cref{lem:skj1-not-boundary}.
\end{proof}

Now, we can establish results similar to \Cref{thm:cv-outer-wk}, \Cref{cor:stationarity}, and \Cref{thm:eps-kkt-stat} for \Cref{alg:bar-outer-z}.

\begin{shadytheorem}
  \label{thm:cv-outer-z}
  Let \Cref{asm:model}, \Cref{asm:liminf-delta-nu} and \Cref{asm:kappa-mdb-eps-cv} be satisfied, and for all \(j\), \((f + h)(x_{k, j}) \ge (f + h)_{\textup{low}, k}\).
  Then, there exists a subsequence \(N \in \mathcal{N}_{\infty}\) such that for all \(k \in N\),
  \begin{equation}
    \label{eq:asympt-stat-zkj}
    -\bar \nu_k^{-1} \bar s_k^{\mathcal{L}} \in \nabla f(\bar x_k) - \bar z_k + \partial \psi(\bar s_k^{\mathcal{L}}; \bar x_k),
  \end{equation}
  with
  \begin{equation}
    \label{eq:nkjk-bar-skjk1-lim}
    \nu_{k,j}^{-1} \|\bar s_k^{\mathcal{L}}\| \le \sqrt{2} \epsilon_{d,k} \rightarrow 0.
  \end{equation}
\end{shadytheorem}

\begin{proof}
  \Cref{lem:approx-stat-kj-z} and \Cref{lem:bar-skj1-not-boundary} lead to~\eqref{eq:asympt-stat-zkj}.
  \Cref{lem:bar-skjk1-epsilondk-ineq} shows that \(\nu_{k,j}^{-1} \|s_{k,j}^{\mathcal{L}}\| \le \sqrt{2}\epsilon_{d,k}\), and, as \(\epsilon_{d,k} \rightarrow 0\),~\eqref{eq:nkjk-bar-skjk1-lim} is satisfied.
\end{proof}

\begin{shadycorollary}
  \label{cor:stationarity-xi-z}
  Under the assumptions of \Cref{thm:cv-outer-z} and \Cref{asm:subgrad-cv}, let \(\bar x\) and \(\bar z\) be limit points of \(\{\bar x_k\}\) and \(\{\bar z_k\}\), respectively.
  Then,
  \begin{equation}
    0 \in \nabla f(\bar x) - \bar z + \partial h(\bar x) \quad \text{and} \quad \bar X \bar Z e = 0.
  \end{equation}
  In this case, when the CQ is satisfied, \(\bar x\) is first-order stationary for~\eqref{eq:nlp}.
\end{shadycorollary}

\begin{proof}
  In \Cref{thm:cv-outer-z}, \(N\) can be chosen such that \(\bar x_k \underset{k \in N}{\longrightarrow} \bar x\) and \(\bar z_k \underset{k \in N}{\longrightarrow} \bar z\).
  We apply \Cref{asm:subgrad-cv} to~\eqref{eq:asympt-stat-zkj} and~\eqref{eq:nkjk-bar-skjk1-lim} to obtain
  \begin{equation*}
    -\nabla f(\bar x) + \bar z \in \partial \psi(0; \bar x) = \partial h(\bar x),
  \end{equation*}
  and we conclude as in the proof of \Cref{cor:stationarity}.
\end{proof}

Finally, we show a result similar to \Cref{thm:eps-kkt-stat} for the \((\epsilon_p, \epsilon_d)\)-KKT optimality.
We use again a subsequence \(N\) as in \cref{thm:cv-outer-z}, \(v_k \in \partial h(\bar x_k + \bar s_k^{\mathcal{L}})\), and we define
\begin{equation}
  \label{eq:def-eps-hk-L}
  \epsilon_{h,k}^{\mathcal{L}} = \dist(-\bar \nu_k^{-1} \bar s_k^{\mathcal{L}} - \nabla f(\bar x_k) + \bar z_k, \, \partial h(\bar x_k + \bar s_{k,1})).
\end{equation}

\begin{shadytheorem}
  \label{thm:eps-kkt-stat-z}
  Let the assumptions of \Cref{thm:cv-outer-z} be satisfied, and \(\epsilon_{h,k}^{\mathcal{L}}\) be defined in~\eqref{eq:def-eps-hk-L}.
  Then, there exists a subsequence \(N \in \mathcal{N}_{\infty}\) such that for all \(k \in N\), \(\bar x_k + \bar s_k^{\mathcal{L}}\) is \((\bar \epsilon_{p,k}^{\mathcal{L}}, \bar \epsilon_{d,k}^{\mathcal{L}})\)-KKT optimal with constants
  \begin{align*}
    \bar \epsilon_{p,k}^{\mathcal{L}} &= \epsilon_{p,k} + \sqrt{n}\mu_k + \sqrt{2}\bar \nu_k\epsilon_{d,k}\|\bar z_k\| \\
    \bar \epsilon_{d,k}^{\mathcal{L}} &= \epsilon_{h,k}^{\mathcal{L}} + \sqrt{2} \epsilon_{d,k}(1 + L_f \bar \nu_k).
  \end{align*}
\end{shadytheorem}

\begin{proof}
  \Cref{thm:cv-outer-z} guarantees that there exists an infinite subsequence \(N\) such that for all \(k \in N\),
  \begin{equation*}
    -\bar \nu_k^{-1} \bar s_k^{\mathcal{L}} - \nabla f(\bar x_k) + \bar z_k \in \partial \psi(\bar s_k^{\mathcal{L}}; \bar x_k),
  \end{equation*}
  Since \(\partial h(\bar x_k + \bar s_k^{\mathcal{L}})\) is closed and nonempty, we can choose \(v_k \in \partial h(\bar x_k + \bar s_k^{\mathcal{L}})\) such that, for \(y_k^{\mathcal{L}} := -\bar \nu_k^{-1} \bar s_k^{\mathcal{L}} - \nabla f(\bar x_k) + \bar z_k\), we have \(\|y_k^{\mathcal{L}}\| = \epsilon_{h,k}^{\mathcal{L}}\).
  Now,
  \begin{equation*}
    (v_k - (-\bar \nu_k^{-1} \bar s_k^{\mathcal{L}} - \nabla f(\bar x_k) + \bar z_k)) - \bar \nu_k^{-1} \bar s_k^{\mathcal{L}} \in \nabla f(\bar x_k) - \bar z_k + \partial h(\bar x_k + \bar s_k^{\mathcal{L}}),
  \end{equation*}
  which can also be written as
  \begin{equation*}
    y_k^{\mathcal{L}} - \bar \nu_k^{-1} \bar s_k^{\mathcal{L}} + \nabla f(\bar x_k + \bar s_k^{\mathcal{L}}) - \nabla f(\bar x_k) \in \nabla f(\bar x_k + \bar s_k^{\mathcal{L}}) - \bar z_k + \partial h(\bar x_k + \bar s_k^{\mathcal{L}}).
  \end{equation*}
  The triangle inequality combined and the Lipschitz constant \(L_f\) of \(\nabla f\) leads to
  \begin{equation*}
    \|y_k^{\mathcal{L}} - \bar \nu_k^{-1} \bar s_k^{\mathcal{L}} + \nabla f(\bar x_k + \bar s_k^{\mathcal{L}}) - \nabla f(\bar x_k)\| \le \epsilon_{h,k}^{\mathcal{L}} + \bar \nu_k^{-1} \|\bar s_k^{\mathcal{L}}\| + L_f \|\bar s_k^{\mathcal{L}}\|.
  \end{equation*}
  \Cref{lem:bar-skjk1-epsilondk-ineq} then implies that
  \begin{equation*}
    \label{eq:eps-d-ineq-kkt-stat-z}
    \|y_k^{\mathcal{L}} - \bar \nu_k^{-1} \bar s_k^{\mathcal{L}} + \nabla f(\bar x_k + \bar s_k^{\mathcal{L}}) - \nabla f(\bar x_k)\| \le \epsilon_{h,k}^{\mathcal{L}} + \sqrt{2} \epsilon_{d,k}(1 + L_f \bar \nu_k).
  \end{equation*}
  Finally, the inequalities of~\eqref{eq:eps-p-ineq-kkt-stat} still hold when replacing \(\bar s_{k,1}\) by \(\bar s_k^{\mathcal{L}}\): for all \(i \in \{1, \ldots, n\}\), \((\bar x_k)_i (\bar z_k)_i \le \epsilon_{p,k} + \sqrt{n} \mu_k + \sqrt{2}\bar \nu_k\epsilon_{d,k}\|\bar z_k\|\).
\end{proof}

When \(\bar x_k \underset{k \in N}{\longrightarrow} \bar x \in \R_{+}^n\) and \(\bar z_k \underset{k \in N}{\longrightarrow} \bar z \in \R_+^n\), \(\{\bar z_k\}_N\) is bounded, thus \(\bar \epsilon_{p,k}^{\mathcal{L}} \underset{N}{\longrightarrow} 0\) in \Cref{thm:eps-kkt-stat-z}.
If \(\epsilon_{h,k}^{\mathcal{L}} \rightarrow 0\), we also have \(\bar \epsilon_{d,k}^{\mathcal{L}} \underset{N}{\longrightarrow} 0\).

\section{Implementation and numerical experiments}%
\label{sec:numerical}

All solvers tested are available from \href{https://github.com/JuliaSmoothOptimizers/RegularizedOptimization.jl}{RegularizedOptimization.jl}.
We define
\begin{equation*}
    \epsilon_{d,k} := \epsilon_k + \epsilon_{r,i} \nu_{k,0}^{-1/2}\xi_{\delta_k}^{\mathcal{L}}(\Delta_{k,0}; x_{k,0}, z_{k,0}, \nu_{k,0})^{1/2} \quad \text{and} \quad \epsilon_{p,k} := \epsilon_k,
\end{equation*}
where \(\epsilon_{r,i} \ge 0\) is a predefined relative tolerance for the inner iterations.
\Cref{alg:bar-inner} terminates when
\begin{equation*}
  \begin{aligned}
  \nu_{k,j}^{-1/2}\xi_{\delta_k}^{\mathcal{L}}(\Delta_{k,j}; x_{k,j}, z_{k,j}, \nu_{k,j_k})^{1/2} & < \epsilon_{d,k} \\
  \|X_{k,j} z_{k,j} - \mu_k e\| & < \epsilon_{p,k}.
  \end{aligned}
\end{equation*}
We use the constant \(\kappa_{\textup{bar}} = 10^6\) for \(\Theta_{k,j}\) in~\eqref{eq:pd-model-theta} and \(\epsilon_{r,i} = 10^{-1}\).
We set \(\epsilon_0 = \mu_0 = 1\), \(\mu_{k+1} = \mu_k / 10\), \(\epsilon_k = \mu_k^{1.01}\), and \(\Delta_{k,0} = 1000\mu_k\), similarly as \citet{conn-gould-orban-toint-2000} did in the smooth case for their interior-point trust-region algorithm.   
In addition, we go to iteration \(k+1\) and set \(x_{k+1,0} = x_{k,j}\) if \Cref{alg:bar-inner} performs more than \(j = 200\) iterations.
Any iteration \(k\) where \(x_{k,0}\) is not first-order stationary for~\eqref{eq:nlp} has \(\epsilon_{d,k} > \epsilon_{p,k} = \epsilon_k\), thus the remarks below \Cref{asm:kappa-mdb-eps-cv} are not sufficient to prove that this assumption is always satisfied, however, we observe satisfying performance with these parameters.
Although it is possible to use \(\epsilon_{d,k} = \epsilon_k\), it would require more inner iterations.

To declare convergence of \Cref{alg:bar-outer-z}, we use the following criteria
\begin{subequations}
  \begin{align}
    \mu_k &< \epsilon \\
    \|\bar X_k \bar z_k - \mu_k e\| &< \epsilon \\
    \bar \nu_k^{-1/2} \xi_{\delta_k}^{\mathcal{L}}(\bar \Delta_k; \bar x_k, \bar z_k, \bar \nu_k)^{1/2} &< \epsilon,
  \end{align}
\end{subequations}
where \(\epsilon = \epsilon_a + \epsilon_r \nu_{1,1}^{-1/2} \xi_{\delta_1}^{\mathcal{L}}(\Delta_{1,1}; x_{1,1}, \nu_{1,1})^{1/2}\), for some \(\epsilon_a \ge 0\).
In our experiments, we chose \(\epsilon_a = \epsilon_r = 10^{-4}\).
We could not base our criteria upon \(\xi^{\mathcal{L}}\) in~\eqref{eq:def-xibar} because we do not know the \(\nu\) associated to this measure (which is different from the \(\nu_{k,j}\) generated by \Cref{alg:bar-inner} associated to the barrier subproblem).

Once \Cref{alg:bar-outer-z} terminates, we use a crossover technique to set \(x_i = 0\) or/and \(z_i = 0\), to respect the complementarity condition.
To do so, we check the final value of \(x_i\) (resp. \(z_i\)), and if it is smaller than \(\sqrt{\mu}\) we set it to zero.
If both \(x_i\) and \(z_i\) are smaller than \(\mu^{1/4}\), we set them to zero.

The subproblems in \Cref{alg:bar-outer-z} are solved with the algorithm R2 \citep{aravkin-baraldi-orban-2022}, and we compare \Cref{alg:bar-outer} to TR \citep{aravkin-baraldi-orban-2022} with R2 used as a subsolver, and R2 used by itself. 
\Cref{alg:bar-outer-z} will be denoted RIPM.
Finally, we introduce a variant of RIPM named RIPMDH (\emph{Regularized Interior Proximal Method with Diagonal Hessian approximations}), that uses the same idea as our algorithm TRDH \citep{leconte-orban-2023}: instead of using LBFGS or LSR1 quasi-Newton approximations for \(B_{k,j}\), we use diagonal quasi-Newton approximations so that~\eqref{eq:tr-sub} with \(\varphi\) defined in~\eqref{eq:pd-model-theta} can be solved analytically for specific seperable regularizers \(h\).
TRDH and RIPMDH use the Spectral Gradient update in all our results.
We choose either \(h(x) = \lambda \|x\|_0\) or \(h(x) = \lambda \|x\|_1\), where \(\lambda > 0\).
When \(h(x) = \lambda \|x\|_0\), RIPM and RIPMDH denote \Cref{alg:bar-outer} instead of \Cref{alg:bar-outer-z} because \(h\) is not convex.

For simplicity, we described how to solve~\eqref{eq:nlp} with the constraint \(x \ge 0\), but RIPM and RIPMDH are actually able to handle box constraints \(\ell \le x \le u\).
These more general constraints can be handled with minor modifications using the barrier function
\begin{equation}
  \label{eq:phik-l-u}
  \tilde \phi_k(x) := -\mu_k \sum_{i=1}^n \log (x_i - \ell_i) - \mu_k \sum_{i=1}^n \log (u_i - x_i)
\end{equation}
instead of \(\phi_k\) \citep[Section~\(13.8\)]{conn-gould-toint-2000}.
When \(\ell_i = -\infty\) (resp. \(u_i = +\infty\)) for some \(i \in \{1, \ldots, n\}\), we remove the term \(\log(x_i - \ell_i)\) (resp. \(\log (u_i - x_i)\)) from the first (resp. second) sum in~\eqref{eq:phik-l-u}.

Our results report
\begin{itemize}
  \item the final \(f(x)\);
  \item the final \(h(x) / \lambda\), where \(\lambda\) is a parameter relative to our regularization function \(h\);
  \item the final stationarity measure \(\sqrt{\xi / \nu}\);
  \item \(\| x - x^*\|\), where \(x_*\) is the exact solution, if it is available;
  \item the number of smooth objective evaluations \(\# f\);
  \item the number of gradient evaluation \(\# \nabla f\);
  \item the number of proximal operator evaluations \(\# \textup{prox}\);
  \item the elapsed time \(t\) in seconds.  
\end{itemize}

Our main goal is to reduce the number of objective and gradient evaluation, as they are typically costly to evaluate.
Since we did not fully optimize the allocations in our algorithms, we do not pay attention to the elapsed time, and we only report it in the tables for information.

Once a problem has been solved by all solvers, we compare their final objective values and we save the smallest, that we denote \((f+h)^*\).
Then, for all solvers, we plot \((f+h)(x_k) - (f+h)^*\) for every iteration \(k\) where the gradient \(\nabla f\) is evaluated.
This allows us to represent the evolution of the objective per gradient evaluation.
For the last gradient evaluation of RIPM and RIPMDH, we display their final objective value after applying the crossover technique.

\subsection{Box-constrained quadratic problem}
\label{sec:box-cstr-qp}

\begin{figure}[ht]%
  \centering
  \includetikzgraphics[width = 0.8\linewidth]{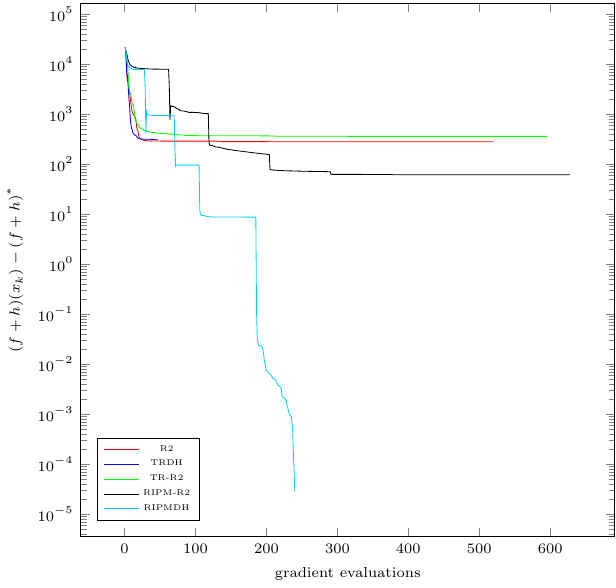}
  \caption{%
    \label{fig:qp-rand-plots}
    Plots of the objective of~\eqref{eq:constrained-qp} per gradient evaluation with different solvers.
  }
\end{figure}

For our first numerical experiment, we solve 

\begin{equation}
  \label{eq:constrained-qp}
  \minimize{x} c^T x + \tfrac{1}{2} x^T H x + h(x) \quad \st \ \ell \le x \le u, 
\end{equation}
which is similar to \citep[Section~\(7.1\)]{shen-xue-zhang-wang-2020}, where \(h = \lambda \|\cdot\|_1\), \(H = A + A^T\), \(A \in \R^{n \times n}\) has nonzero components with probability \(p = 10^{-4}\) following a normal law of mean \(0\) and standard deviation \(1\), \(c \in \R^n\) has components generated using a normal distribution of mean \(0\) and standard deviation \(1\), \(\ell = -e - t_{\ell}\) and \(u = e + t_u\), with \(t_{\ell} \in \R^n\), \(t_u \in \R^n\) are vectors sampled from a uniform distribution between \(0\) and \(1\).
We chose \(n = 10^5\), and use a LSR1 quasi-Newton approximation for TR and RIPM.
For \(\lambda \ge 1.0\), the components \(x_i\) of the solutions returned by TR, TRDH, R2, RIPM and RIPMDH satisfy \(x_i \in \{\ell_i, u_i, 0\}\) for almost all \(i \in \{1, \ldots, n\}\).
In this case, we observe that TR, TRDH and R2 are more efficient than RIPM.
However, as we decrease \(\lambda\), we get more components \(x_i \notin \{\ell_i, u_i, 0\}\).
We show results with \(\lambda = 10^{-1}\) in \Cref{fig:qp-rand-plots} and \Cref{tbl:qp-rand}.
TRDH performs the least amount of objective, gradient and proximal operators evaluations.
RIPMDH finds the smallest final objective value, and performs fewer objective, gradient and proximal operator evaluations than TR-R2.
RIPM terminates with a criticality measure higher than the other solvers, but we observe that its final objective is smaller than those of TR, TRDH and R2.
For RIPM and RIPMDH, we can clearly see plateaus that delimit the outer iterations.
RIPM-R2 performs many more proximal operator evaluations than RIPMDH, because it uses up to \(200\) R2 iterations to solve~\eqref{eq:tr-sub}.
The number of proximal operator evaluations with RIPM-R2 is also much higher than that of TR-R2, because the subproblems solved with R2 in RIPM-R2 have their objective based upon~\eqref{eq:pd-model-theta}, which is not well conditioned when some components of \(x_{k,j}\) approach \(0\), whereas the subproblems in TR-R2 are based upon~\eqref{eq:varphi-f}.

\begin{table}[ht]%
  \centering
  \scriptsize
  \caption{%
  \label{tbl:qp-rand}
  Statistics of~\eqref{eq:constrained-qp}.
  TR and RIPM use an LSR1 Hessian approximation.
  The maximum number of objective evaluations is set to \(800\).
  }
  \begin{tabular}{rrrrrrrr}
    \hline\hline
    \textbf{solver} & \textbf{$f(x)$} & \textbf{$h(x) / \lambda$} & \textbf{$\sqrt{\xi / \nu}$} & \textbf{$\# f$} & \textbf{$\# \nabla f$} & \textbf{$\# prox$} & \textbf{$t$ ($s$)} \\\hline
    R2 & \(-2.29\)e\(+04\) & \(1.5\)e\(+04\) & \(8.5\)e\(-03\) & \(679\) & \(520\) & \(679\) & \(6.8\)e\(-01\) \\
    TRDH & \(-2.28\)e\(+04\) & \(1.5\)e\(+04\) & \(5.9\)e\(-05\) & \(57\) & \(47\) & \(113\) & \(3.6\)e\(-01\) \\
    TR-R2 & \(-2.28\)e\(+04\) & \(1.5\)e\(+04\) & \(9.9\)e\(-03\) & \(801\) & \(596\) & \(12639\) & \(8.5\)e\(+00\) \\
    RIPM-R2 & \(-2.30\)e\(+04\) & \(1.4\)e\(+04\) & \(3.4\)e\(+00\) & \(801\) & \(628\) & \(101019\) & \(5.0\)e\(+01\) \\
    RIPMDH & \(-2.32\)e\(+04\) & \(1.5\)e\(+04\) & \(8.7\)e\(-03\) & \(313\) & \(241\) & \(628\) & \(2.7\)e\(+00\) \\\hline\hline
  \end{tabular}
\end{table}

\subsection{Sparse nonnegative matrix factorization (NNMF)}
\label{sec:nnmf}

The second experiment considered is the sparse nonnegative matrix factorization (NNMF) problem from \citet{kim-park-2008}.
Let \(A \in \R^{m \times n}\) have nonnegative entries.
Each column of \(A\) represents an observation, and is generated using a mixture of Gaussians where negative entries are set to zero.
We factorize \(A \approx WH\) by separating \(A\) into \(k < \min(m, n)\) clusters, where \(W \in \R^{m \times k}\), \(H \in \R^{k \times n}\) both have nonnegative entries and \(H\) is sparse.
This problem can be written as
\begin{equation}%
  \label{eq:nnmf}
  \minimize{W, H} \tfrac{1}{2} \|A - W H\|_F^2 + h(H) \quad \st \ W, H \ge 0,
\end{equation}
where \(h(H) = \lambda \|\textup{vec}(H)\|_1\) and \(\textup{vec}(H)\) stacks the columns of \(H\) to form a vector.

\begin{figure}[ht]%
  \centering
  \includetikzgraphics[width = 0.8\linewidth]{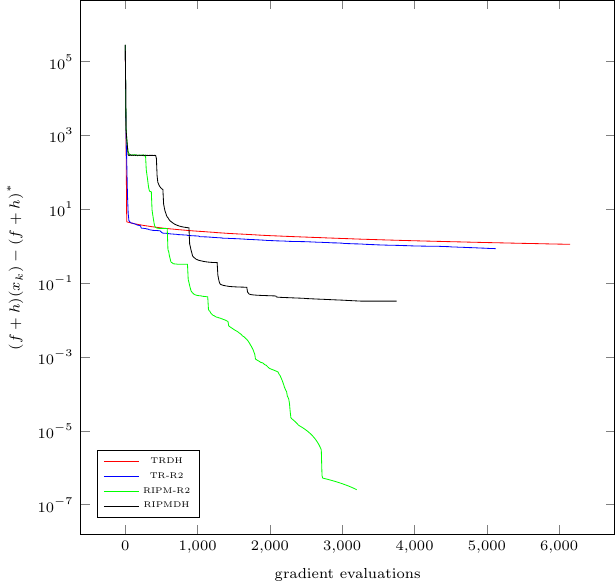}
  \caption{%
    \label{fig:nnmf-plots}
    Plots of the objective of~\eqref{eq:nnmf} per gradient evaluation with different solvers.
  }
\end{figure}

We set \(m = 100\), \(n = 50\), \(k = 5\), \(\lambda = 10^{-1}\), and report the statistics in \Cref{tbl:nnmf}.
For this particular problem, we use \(\epsilon_r = 10^{-6}\), which allows for more accurate solves and for a better visualization of the evolution of the objective values, shown in \Cref{fig:nnmf-plots}.
We observe that RIPM-R2 and RIPMDH are the only solvers to terminate.
They outperform R2 and TR-R2 in terms of number of objective and gradient evaluations, and their final objective value is also smaller.
R2 and TR-R2 reach the maximum number of iterations.
The objective of RIPM-R2 and RIPMDH is higher than that of R2 and TR-R2 only in the early iterations, because the barrier function has more effect when \(\mu_k\) is larger.
RIPM-R2 performs less objective and gradient evaluations than RIPMDH, but much more proximal operator evaluations because of the reasons evoked in~\Cref{sec:box-cstr-qp}.

\begin{table}[ht]
  \centering
  \scriptsize
  \caption{%
  \label{tbl:nnmf}
  Statistics of~\eqref{eq:nnmf}.
  TR and RIPM use an LSR1 Hessian approximation.
  The maximum number of objective evaluations is set to \(8000\).}
  \begin{tabular}{rrrrrrrr}
    \hline\hline
    \textbf{solver} & \textbf{$f(x)$} & \textbf{$h(x) / \lambda$} & \textbf{$\sqrt{\xi / \nu}$} & \textbf{$\# f$} & \textbf{$\# \nabla f$} & \textbf{$\# prox$} & \textbf{$t$ ($s$)} \\\hline
    TRDH & \(1.25\)e\(+02\) & \(3.1\)e\(+01\) & \(1.6\)e\(-01\) & \(8001\) & \(6156\) & \(16000\) & \(2.0\)e\(+00\) \\
    TR-R2 & \(1.25\)e\(+02\) & \(2.8\)e\(+01\) & \(1.2\)e\(-01\) & \(8001\) & \(5122\) & \(150563\) & \(6.4\)e\(+00\) \\
    RIPM-R2 & \(1.25\)e\(+02\) & \(1.9\)e\(+01\) & \(2.5\)e\(-02\) & \(4501\) & \(3210\) & \(470975\) & \(1.1\)e\(+01\) \\
    RIPMDH & \(1.25\)e\(+02\) & \(2.0\)e\(+01\) & \(1.4\)e\(-02\) & \(4602\) & \(3759\) & \(9205\) & \(1.3\)e\(+00\) \\\hline\hline
  \end{tabular}
\end{table}

\subsection{FitzHugh-Nagumo problem (FH)}
\label{sec:fh}

\begin{figure}[ht]%
  \centering
  \includetikzgraphics[width = 0.8\linewidth]{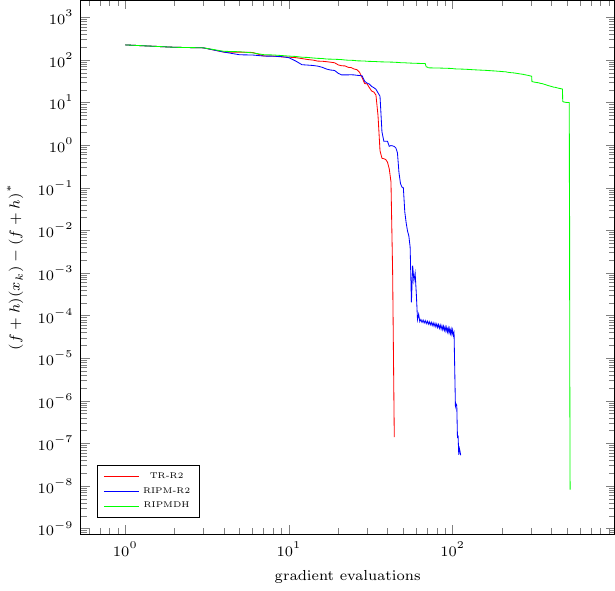}
  \caption{%
    \label{fig:fh-plots}
    Plots of the objective of~\eqref{eq:fh-cstr} per gradient evaluation with different solvers.
  }
\end{figure}

We sample the functions \(V(t; x)\) and \(W(t; x)\) satisfying the \citet{fitzhugh-1995} and \citet{nagumo-arimoto-1962} model for neuron activation, where \(x \in \R^5\), as \(v(x) = (v_1(x), \ldots, v_{n+1}(x))\) and \(w(x) = (w_1(x), \ldots, w_{n+1}(x))\).   
\begin{equation}
  \label{eq:fh-nagumo}
  \frac{\mathrm{d} V}{\mathrm{d} t} = (V - V^3/3 - W + x_1) x_2^{-1}, \quad \frac{\mathrm{d} W}{\mathrm{d} t} = x_2 (x_3 V  - x_4 W + x_5).
\end{equation}
The time interval \(t \in [0, 20]\) is discretized, with the initial conditions \((V(0), W(0)) = (2, 0)\).
We solve
\begin{equation}%
  \label{eq:fh-cstr}
  \minimize{x} \tfrac{1}{2} \| (v(x) - \bar v (\bar x), w(x) - \bar w (\bar x) )\|_2^2 + h(x), \quad \st \ x_2 \ge 0.5,
\end{equation}
where \(h(x) = \lambda \|x\|_0\) with \(\lambda = 10\), \(n = 100\), and report the statistics in \Cref{tbl:fh-cstr}.
Since \(h\) is not convex, we use \Cref{alg:bar-outer} instead of \Cref{alg:bar-outer-z}.
We do not show results with R2 because it encounters a numerical error during the solve of a differential equation to compute the objective. 
The evolution of the objective per gradient evaluation is shown in \Cref{fig:fh-plots}.
To improve readability, we choose to show the number of gradient evaluations on a logarithmic scale, and not to plot results with TRDH.
All solvers converge to the value \((0.00, 0.50, 0.54, 0.00, 0.00)\) except for TRDH that has a higher final objective value than the other solvers.
TR-R2 is the fastest, and seems the most suited to solve smaller problems such as~\eqref{eq:fh-cstr}.
RIPM and RIPMDH still converge, but the latter is much slower.
However, RIPMDH performs the least amount of proximal operator evaluations.

\begin{table}[ht]
  \centering
  \scriptsize
  \caption{%
  \label{tbl:fh-cstr}
  Statistics of~\eqref{eq:fh-cstr}.
  TR and RIPM use an LBFGS Hessian approximation.
  }
  \begin{tabular}{rrrrrrrr}
    \hline\hline
    \textbf{solver} & \textbf{$f(x)$} & \textbf{$h(x)/\lambda$} & \textbf{$\sqrt{\xi / \nu}$} & \textbf{$\# f$} & \textbf{$\# \nabla f$} & \textbf{$\# prox$} & \textbf{$t$ ($s$)} \\\hline
    TRDH & \(6.05\)e\(+00\) & \(3\) & \(2.9\)e\(+01\) & \(1001\) & \(697\) & \(2000\) & \(4.5\)e\(+00\) \\
    TR-R2 & \(4.40\)e\(+00\) & \(2\) & \(4.8\)e\(-03\) & \(53\) & \(45\) & \(4627\) & \(3.3\)e\(-01\) \\
    RIPM-R2 & \(4.40\)e\(+00\) & \(2\) & \(1.3\)e\(-02\) & \(261\) & \(112\) & \(10139\) & \(1.0\)e\(+00\) \\
    RIPMDH & \(4.40\)e\(+00\) & \(2\) & \(1.5\)e\(-02\) & \(798\) & \(523\) & \(1600\) & \(3.8\)e\(+00\) \\\hline\hline
  \end{tabular}
\end{table}

\subsection{Constrained basis pursuit denoise (BPDN)}
\label{sec:bpdn-cstr}

We solve the basis pursuit denoise problem (BPDN) \citep{tibshirani-1996,donoho-2006} with additional bound constraints.
Let \(m = 200\), \(n = 512\), \(b = A x_\star + \epsilon\), where \(\epsilon \sim \mathcal{N}(0, 0.01)\), \(A \in \R^{m \times n}\) has orthonormal rows, and \(x_\star\) is a vector of zeros, except for \(5\) of its components that are set to \(1\).
The constrained BPDN problems is written as
\begin{equation}
  \label{eq:bpdn-cstr}
  \minimize{x} \tfrac{1}{2} \|A x - b\|_2^2 + h(x) \quad \st \ x \ge 0,
\end{equation}
where \(h(x) = \lambda \|x\|_1\).
We use \(\lambda = \|A^T b\|_{\infty} / 10\).

\begin{figure}[ht]%
  \centering
  \includetikzgraphics[width = 0.8\linewidth]{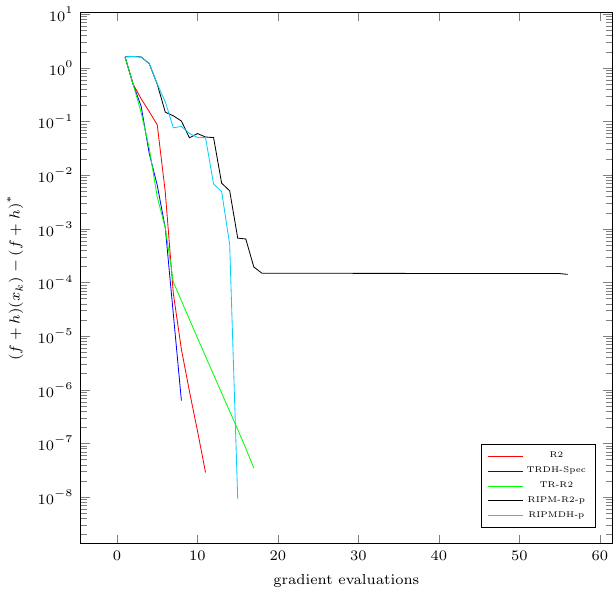}
  \caption{%
    \label{fig:bpdn-plots}
    Plots of the objective of~\eqref{eq:bpdn-cstr} per gradient evaluation with different solvers.
  }
\end{figure}

The statistics are shown in \Cref{tbl:bpdn-cstr}.
R2, TRDH and TR-R2 are much more efficient than RIPM on this problem.
This could come from the fact that there are many active bounds in the solution.
However, this was also the case for the NNMF problem of \Cref{sec:nnmf}, for which RIPM seems more efficient.
Further investigations should seek to understand such behaviours on different problems.
RIPM-R2-p and RIPMDH-p use the modifications \(\mu_0 = 10^{-3}\) and \(\epsilon_{r,i} = 1.0\), which make RIPMDH on~\eqref{eq:bpdn-cstr} surpass TR-R2 and close to TRDH.
\Cref{fig:bpdn-plots} shows the evolution of the objective values.
RIPM and RIPMDH are not included to improve readability.
\begin{table}[ht]
  \centering
  \scriptsize
  \caption{%
  \label{tbl:bpdn-cstr}
  Statistics of~\eqref{eq:bpdn-cstr}.
  TR and RIPM use an LSR1 Hessian approximation.
  }
  \begin{tabular}{rrrrrrrrr}
    \hline\hline
    \textbf{solver} & \textbf{$f(x)$} & \textbf{$h(x)/\lambda$} & \textbf{$\sqrt{\xi / \nu}$} & \textbf{$\|x-x^*\|_2$} & \textbf{$\# f$} & \textbf{$\# \nabla f$} & \textbf{$\# prox$} & \textbf{$t$ ($s$)} \\\hline
    R2 & \(3.91\)e\(-02\) & \(8.7\)e\(+00\) & \(1.8\)e\(-04\) & \(4.1\)e\(-01\) & \(11\) & \(11\) & \(11\) & \(4.0\)e\(-03\) \\
    TRDH & \(3.91\)e\(-02\) & \(8.7\)e\(+00\) & \(6.8\)e\(-05\) & \(4.1\)e\(-01\) & \(9\) & \(9\) & \(17\) & \(8.0\)e\(-03\) \\
    TR-R2 & \(3.91\)e\(-02\) & \(8.7\)e\(+00\) & \(1.6\)e\(-04\) & \(4.1\)e\(-01\) & \(17\) & \(17\) & \(35\) & \(1.1\)e\(-02\) \\
    RIPM-R2 & \(4.12\)e\(-02\) & \(8.7\)e\(+00\) & \(7.1\)e\(-03\) & \(4.2\)e\(-01\) & \(915\) & \(915\) & \(12143\) & \(3.3\)e\(+00\) \\
    RIPMDH & \(3.92\)e\(-02\) & \(8.7\)e\(+00\) & \(9.9\)e\(-04\) & \(4.1\)e\(-01\) & \(361\) & \(225\) & \(724\) & \(2.1\)e\(-01\) \\
    RIPM-R2-p & \(4.22\)e\(-02\) & \(8.7\)e\(+00\) & \(1.0\)e\(-02\) & \(4.3\)e\(-01\) & \(56\) & \(56\) & \(8331\) & \(1.5\)e\(-01\) \\
    RIPMDH-p & \(3.91\)e\(-02\) & \(8.7\)e\(+00\) & \(5.6\)e\(-04\) & \(4.1\)e\(-01\) & \(15\) & \(15\) & \(32\) & \(9.0\)e\(-03\) \\\hline\hline
  \end{tabular}
\end{table}


  \section{Discussion and future work}

We have presented RIPM, a trust-region interior-point method to solve nonsmooth regularized problems with box constraints, and RIPMDH, a variant based upon techniques of \citet{leconte-orban-2023}.
These algorithms solve a sequence of unconstrained barrier subproblems to obtain a sequence of approximate solutions of~\eqref{eq:nlp}.
We have shown the convergence of the inner barrier subproblems, and we have characterized the degree of \((\epsilon_p, \epsilon_d)\)-KKT optimality for every outer iteration past a certain rank.
Under the assumption that the iterates remain bounded, we have shown that RIPM converges to a first-order stationary point for~\eqref{eq:nlp}.
We compared RIPM and RIPMDH to projected-direction methods with a separable regularization function.

RIPM and RIPMDH perform well on the box-constrained quadratic problem of \Cref{sec:box-cstr-qp} and on the NNMF problem of \Cref{sec:nnmf}.
They are not as efficient on the FH problem of \Cref{sec:fh} and the constrained BPDN problem of \Cref{sec:bpdn-cstr}, which may suggest that projected-direction methods may be more efficient to solve problems with fewer variables and constraints.
However, as observed with RIPM-R2-p and RIPMDH-p, the modification of two parameters of RIPM and RIPMDH improves their efficiency significantly on the constrained BPDN problem.
This suggests that our implementation could benefit from parameter tuning.

Future work may include generalizing the algorithm to constraints of the form \(c_i(x) \le 0\) with \(i \in \{1, \ldots, m\}\) for some \(m > 0\), where the \(c_i\) are continuously differentiable and Lipschitz-gradient continuous, as in \citep[Section~\(13.9\)]{conn-gould-toint-2000} in the smooth case, or \citep{demarchi-themelis-2022} for nonsmooth problems.

Another improvement would be to scale the trust region to allow greater search directions along the boundary of the feasible domain.
This is explained more in detail in \citep[Section~\(13.7\)]{conn-gould-toint-2000} for trust-regions based upon the \(\ell_2\)-norm.
However, we could not find an alternative for trust-regions based upon the \(\ell_{\infty}\)-norm that led to satisfying numerical results.

In \Cref{sec:new-crit-meas}, \Cref{asm:psi-convex} does not allow the use of \(h = \|\cdot\|_0\) with \Cref{alg:bar-outer-z}.
It would be interesting to see whether it is possible to establish convergence properties similar to those of \Cref{alg:bar-outer} without this assumption.
One way to do this might be to replace \(\xi_{\textup{cp}}\) by \(\xi_{\delta_k}^{\mathcal{L}}\) in~\eqref{eq:cauchy-decrease} when \(j\) is large enough, but we did not manage to justify that this change results in a reasonable \Cref{asm:cauchy-decrease}. 

The extension of the convergence results to locally Lipschitz-gradient continuous functions \(f\) could also be studied, based upon the work of \citep{demarchi-themelis-2022b,demarchi-themelis-2022,kanzow-mehlitz-2022}.

Finally, when \(\{B_{k,j}\}_j\) grows unbounded, we may still be able to prove the convergence of RIPM using the analysis of \citet{leconte-orban-2023b}, provided that the norm of the Hessian approximations do not grow too fast.

  \bibliographystyle{abbrvnat}
  \bibliography{abbrv,references}


  

\end{document}